\documentclass[11pt]{article}
\usepackage{epsfig}
\usepackage{amsmath}
\usepackage{amsfonts}
\usepackage{mathrsfs}
\usepackage{amsmath, amssymb}
\usepackage{graphicx}
\usepackage{geometry}

 \def\cal#1{\mathcal{#1}}

\makeatletter
\@tfor\next:=abcdefghijklmnopqrstuvwxyzABCDEFGHIJKLMNOPQRSTUVWXYZ\do{%
  \def\command@factory#1{%
    \expandafter\def\csname cal#1\endcsname{\mathcal{#1}}
    \expandafter\def\csname frak#1\endcsname{\mathfrak{#1}}
    \expandafter\def\csname scr#1\endcsname{\mathscr{#1}}
    \expandafter\def\csname bb#1\endcsname{\mathbb{#1}}
  }
 \expandafter\command@factory\next
}
\makeatother

\newtheorem{thm}{Theorem}

\newtheorem{lem}[thm]{Lemma}
\newtheorem{cor}[thm]{Corollary}
\newtheorem{pro}[thm]{Proposition}

\newtheorem{defn}[thm]{Definition}

\newtheorem{rmk}[thm]{Remark}

\def\square{\hfill${\vcenter{\vbox{\hrule height.4pt \hbox{\vrule width.4pt
height7pt \kern7pt \vrule width.4pt} \hrule height.4pt}}}$}

\newenvironment{pf}{{\it Proof:}\quad}{\square \vskip 12pt}

\title{\bf Relative hyperbolicity for automorphisms of  free products, and free groups}
\author{Fran\c{c}ois Dahmani, and Ruoyu Li}
\date{}

\begin{document}
\maketitle

\begin{abstract} We prove that for a free product $G$ with free factor system $\calG$, any automorphism $\phi$ preserving $\calG$,  atoroidal (in a sense relative to $\calG$) and none of whose power send two different conjugates of subgroups in $\calG$ on conjugates of themselves by the same element,  gives rise to a semidirect product $G\rtimes_\phi \mathbb{Z}$  that is relatively hyperbolic with respect to suspensions of groups in $\calG$. We recover a theorem of Gautero-Lustig and Ghosh that, if $G$ is a free group, $\phi$ an automorphism of $G$, and $\calG$ is its family of polynomially growing subgroups, then the semidirect product by $\phi$ is relatively hyperbolic with respect to the suspensions of these subgroups. We apply the first result to the conjugacy problem for certain   automorphisms (atoroidal and toral)  of free products of abelian groups.   
  \end{abstract}

{\footnotesize
\setcounter{tocdepth}{2}
\tableofcontents

}

\section*{Introduction}

\subsection*{Context}

Given a group $G$ and an automorphism $\phi$, the geometry of $\phi$ is encoded in the group $G\rtimes_\phi \mathbb{Z}$. An appealing feature of geometry in groups  is hyperbolicity, in the sense of Gromov. In some circumstances, it can be exhibited by $\phi$.  One of the most iconic examples of this comes from diffeomeorphisms of surfaces.  Let $\Sigma$ be a closed orientable surface of genus larger than $2$, and $f:\Sigma \to \Sigma$ a pseudo-Anosov diffeomorphism, fixing a base point $x$. Instead of introducing the language of laminations, let us settle for Thurston's characterisation of pseudo-Anosov diffeomorphisms, as those that fail to preserve any finite collection of homotopy classes of simple closed curves.      Thurston \cite{Thur} famously proved  that the mapping torus $M$ of $\Sigma$ by $f$ is a closed three-manifold admitting a hyperbolic metric.     Its fundamental group is therefore a uniform lattice in $PSL_2(\mathbb{C})$ and  it is Gromov-hyperbolic.       Let $\phi$ be the automorphism of the fundamental group $\pi_1(\Sigma, x)$  induced by $f$.  The  fundamental group of $M$ can be expressed as  $\pi_1(M, x) \simeq \pi_1(\Sigma, x) \rtimes_\phi \mathbb{Z}$. Rephrasing this with eyes only for this semi-direct product, we can say that 
  $\pi_1(\Sigma, x) \rtimes_\phi \mathbb{Z}$ is  a hyperbolic group if and only if  $\phi$ is an atoroidal automorphism, in the sense that neither it, nor any of its proper powers, preserve a non-trivial conjugacy class.

In the case of a free group $F$, Brinkmann proved an analoguous result:  $F \rtimes_\phi \mathbb{Z}$ is hyperbolic if and only if $\phi$ is atoroidal. In \cite{DTAMS} the first named author illustrated an application of this geometric study to the conjugacy problem for these automorphisms. 

In this paper, we investigate the case of automorphisms of a free product of groups, $G= H_1*\dots H_p*F_k$.  The conjugacy classes of the groups $H_i$ form a free factor system, and the group $F_k$ is free of rank $k$.  This context is a rich source of examples. It can be traced as early as the work of Fouxe-Rabinovich \cite{FR} (and earlier Golowin and Szadowsky \cite{GSz} for the case of few factors), and it received a modern, illuminating reference in the work of Guirardel and Levitt \cite{GL_out}.       We formulate a relevant property of atoroidality for automorphisms preserving a free factor system. We investigate the geometry of the semi-direct products produced by atoroidal automorphisms.  One cannot expect hyperbolicity in general, due to the nature of the factors of the free product. Thus we are interested in relative hyperbolicity. Interestingly, one cannot expect  proper relative hyperbolicity  in general neither.

Let us illustrate this. Consider a free product $G= \mathbb{Z}^2*\mathbb{Z}^3$, and $\phi$ and automorphism. It has to preserve the conjugacy class of  the factor $\mathbb{Z}^3$, so we may just assume, by conjugating it (by moving the base point)  that it preserves the group  $\mathbb{Z}^3$ seen as a subgroup of $G$. Being an automorphism, it has to send the factor  $\mathbb{Z}^2$ on a conjugate of it by an element of  $\mathbb{Z}^3$. So after conjugating again, we may assume that $\phi$ preserves both $\mathbb{Z}^2$ and $\mathbb{Z}^3$. On each of them, it induces an automorphism described by a matrix in $GL_2(\mathbb{Z})$ and  $GL_3(\mathbb{Z})$. One can choose these matrices so that $\phi$ is genuinely atoroidal. However,   $G\rtimes_\phi \mathbb{Z}$ cannot be interestingly relatively hyperbolic: if it were, the factor  $\mathbb{Z}^3$ would have to be parabolic, hence the semidirect factor $\mathbb{Z}$, normalizing it,  would be in the same parabolic subgroup. But it also normalises the factor $\mathbb{Z}^2$, so this factor would be in the same parabolic subgroup too, and we see that $G\rtimes_\phi \mathbb{Z}$ would be a single parabolic subgroup.  

\subsection*{Results}

We propose conditions for the relative hyperbolicity of  $G\rtimes_\phi \mathbb{Z}$  to hold, when $G= H_1*\dots *H_p*F_k$. Francaviglia and Martino have defined the notion of full irreducibility for automorphisms of free products \cite{SFAM}.  A first condition is the full irreducibility  with the atoroidality, when the Scott complexity  $ (k,p)$  is sufficiently large.  
 We propose the obstruction of twinned pair of subgroups (which is crucial for the case of reducible automorphisms). See Definition \ref{def;twinned}, that we anticipate here: we say that $\phi$ has a twinned pair of subgroups in the free factor system of $G$ if there are two conjugates of free factors $A,B$, and an element $g$, such that, for some $n$,  ${\rm ad}_g \circ \phi^n$    simultaneously preserves $A$ and $B$ (here ${\rm ad}_g$ is the conjugation by $g$).

\begin{thm}\label{thm;rel-hyp-fully-irreducible_intro} (See Theorem \ref{thm;rel-hyp-fully-irreducible})

  Let $G$ be a finitely generated group with a free factor system $\calG$, of Scott complexity $(k,p)$, different from $(1,1)$ and $(0,2)$.      
   Let $\phi\in {\rm  Aut}(G,\calG)$ be  fully irreducible and atoroidal. Assume that it has no twinned subgroups in $\calG$ for $\phi$. 

 Then the semi-direct product $G\rtimes_{\phi}\mathbb{Z}$ is relatively hyperbolic, with respect to the mapping torus of $\calG$.     
\end{thm}

The assumption on the Scott complexity is ensured if  the free product  decomposition of $G$ is different from a single free product $H_1*H_2$, or an HNN extension $H_1*_{\{1\}} = H_1* F_1$. We think, but did not prove, that  the assumption on absence of twinned subgroups always hold for fully irreducible atoroidal automorphisms in the Scott complexity of the statement. In principle though, the group generated by two preserved subgroups conjugated to free factors can fail to be a free factor.

In the reducible case, we prove the following, in which the assumption on absence of twinned subgroups is crucial.

\begin{thm}\label{thm;no_twin_intro} (See Theorem \ref{thm;no_twin})

  Let $G$ be a finitely generated group, and $\calG$ be a free factor system. Let $\phi \in {\rm Aut} (G,\calG)$ be atoroidal for $\calG$.  Assume that there is no pair of twinned subgroups in $G$ for $\calG$ and $\phi$.

Then $G\rtimes_\phi \mathbb{Z}$ is relatively hyperbolic with respect to the mapping torus of $\calG$. 
\end{thm}

This last statement contains Brinkmann's result on free groups. It has a number of potentially interesting cases to which it can be applied, as for instance, the following corollary in which we say that an automorphism $\phi\in {\rm Aut}(G,\calG)$ is toral if  for each $H$ such that $[H]\in \calG$ there exists $g\in G$ such that ${\rm ad}_g\circ \phi|_H$ is the identity on $H$.

\begin{cor}\label{coro;toral_intro} (See Coro. \ref{coro;toral})
  
Assume that $G$ is finitely generated, and that $\calG$ is a free factor system of $G$ consisting of torsion free abelian groups.
If $\phi$ is atoroidal, and toral, then the group $G\rtimes_\phi \mathbb{Z}$ is toral relatively hyperbolic.  
\end{cor}

Theorem \ref{thm;no_twin_intro} can also be applied in the realm of free groups, to  automorphisms of free groups that are non-necessarily atoroidal, by hiding the lack of atoroidality in some subgroups, the so-called polynomially growing subgroups for the automorphism. In the favorable case that these subgroups form a free factor system, one may apply our theorem. In the  general case, we apply a variant of our study, that applies, not only to trees, but to hyperbolic coned-off graphs from these trees. 
We thus  recover the following theorem of Gautero and Lustig, and Ghosh.  See Section \ref{sec;321} for a definition and references about polynomially growing subgroups.

\begin{thm}(Gautero-Lustig, and Ghosh)\label{thm;GLG_intro} (See Theorem \ref{thm;GLG})

  If $\phi$ is an automorphism of a finitely generated free group $F$.    
The semidirect product $F\rtimes_\phi \mathbb{Z}$ is relatively hyperbolic with respect to the mapping torus of the collection of maximal polynomially growing subgroups for the outer class of $\phi$. %
\end{thm}  

In the case where $\phi$ has at least one exponentially growing element (the only interesting case), this theorem   is important. It allows a  natural approach to several problems, otherwise rather complicated (\cite{BriGro} \cite{CouHil} \cite{Djap}  among other).  
Unfortunately, the proof of Gautero and Lustig has to be qualified as incomplete for the time being (it relies on some unwritten, or unavailable material about a certain type of train tracks). Very recently    Ghosh proposed a proof \cite{Gho}, relying on some other, advanced, type of train tracks for free group automorphisms. We notice that the train track technology involved in our proof is significantly more accessible, as it is fundamentally based on Bestvina's elegant approach in \cite{BestBers}, chosen by Francaviglia and Martino \cite{SFAM} for their adaptation for free products.  

As an application of this work, we propose a solution to the conjugacy problem for toral atoroidal automorphisms of free products of non-cyclic free abelian groups, for the Scott complexity  $(0,p)$.

\begin{thm} \label{thm;CP_toral_intro} (See Theorem \ref{thm;CP_toral})

  Let $G$ be a finitely generated free product of non-cyclic free abelian groups, $G=A_1*\dots*A_p$. Denote by $\calA$ the free factor system $\{[A_i], i=1,\dots p\}$.  
  There is an algorithm that,   given $\phi_1, \phi_2$, two automorphisms of $(G, \calA)$ that are atoroidal, and toral, determines whether they are conjugate in $Out(G,\calA)$. 
\end{thm}

\subsection*{On proofs}

This strategy for Theorem \ref{thm;rel-hyp-fully-irreducible_intro}  is to realise the automorphism of the group as an equivariant continuous map from a tree to itself, with special properties. These special properties are those of  a train track map, i.e. for which the cancellations of the paths that are image of two consecutive edges is rather well controlled.   Being a train track map from a tree to itself  allows to iterate the map without losing much in terms of cancellation in reduction of paths.      Train track maps are an invention of Bestvina  and Handel, for free groups, that replaces, in a much non-linear setting, the property of being in Jordan form for a matrix.  A beautiful construction of such maps, by Bestvina \cite{BestBers}, was   adapted by Francaviglia and Martino \cite{SFAM}  to the case of fully irreducible automorphisms of free products of groups.   Using the existence of such a map, one can follow arguments of Brinkmann, and Bestvina Feighn and Handel to prove that, for     $\phi$ an automorphism that is fully irreducible, and that has a certain atoroidality property,    iterating  $\phi$ or $\phi^{-1}$ on hyperbolic elements make their length grow exponentially, and iterating $\phi$ or $\phi^{-1}$ on different conjugates of the free factors make their distance, in a certain sense,   grow exponentially. 
This is sufficient to apply a combination theorem, as in \cite{MjR}, or \cite{GW}, that ensure relative hyperbolicity.

For Theorem \ref{thm;no_twin_intro}, we use an induction on Scott complexity, proving first a relative hyperbolicity for a larger free factor system, for which a power of $\phi$ is fully irreducible, and telescoping with a relative hyperbolic structures for the large free factors of this system. In case of sufficiently large complexity, one can use the previous theorem. In case of low complexity, one proves by hand the relative hyperbolicity, using the combination theorem of \cite{Dah03}. We then prove that a certain condition, the central condition, satisfied by toral automorphisms, ensures that there is no twinned subgroups by $\phi$, our only obstruction.

For  Theorem \ref{thm;GLG_intro}, the argument is similar, but applied, not on the Bass-Serre tree of some decomposition, but on its cone-off over the family of maximal polynomially growing subgroups. Atoroidality, and absence of twinned subgroups in that cone-off space are automatic. We need to introduce the theory of polynomially growing subgroups  (under an automorphism) in the context of a free product, in order to pursue the general strategy that consists in treating fully irreducible automorphisms, and concluding by telescoping the relative hyperbolic structures.

A comment is in order at this stage, even though it is independent of the rest of the paper.  Vincent Guirardel has informed us that a sensible use of the Rips machine on the limit $\bbR$-tree of an automorphism reveals the following. If an automorphism of a free product $G$ is fully irreducible, then either $G$ is a surface group and its polynomially growing subgroups correspond to the boundary components of the surface, or its polynomially growing subgroups are conjugate of the free factors. (In order to see this, we refer to Horbez' \cite[\S 4.2]{TitsAlt}: if the limit tree is so-called relatively free, by definition of relative freeness in  \cite[\S 1.3]{TitsAlt},  all polynomially growing subgroups are conjugate to subgroups of the free factors, and if it is not relatively free, by full irreducibility and  \cite[Prop. 4.11]{TitsAlt}, and finiteness of the number of orbits of points with non-trivial stabilizer, the alternative of \cite[Lemma 4.6]{TitsAlt} reduces to the case of an arational surface tree).
 Thus, in that case, the cone-off of the tree is trivial (one cones-off single vertices).  This insight can be used to remove all mention of coning-off the subtrees  in the proof of Theorem \ref{thm;GLG_intro}, and use only in the Bass-Serre trees of the free product.

Theorem \ref{thm;CP_toral_intro} largely follows the approach of \cite{DTAMS}. 
Given the work in \cite{DT}, it is tempting, and might be possible, to try to extend this result to the case where $A_1, \dots, A_p$ are nilpotent and $\phi_i$ induce the identity on them, or  to  the case where the $A_i$ are abelian and $\phi_i$ induce unipotent automorphism on them (after suitable conjugation). There are some difficulties though, as the current lack of computability of the automorphism group of the semi-direct products.

\subsection*{Acknowledgments}

We wish to thank  I. Chatterji, S. Francaviglia,   F. Gautero,  V. Guirardel, M. Mj, and A. Parreau for discussions related to this work. 

Ruoyu Li was supported by the Chinese Scholarship Council.

\numberwithin{thm}{section}

\section{Setting and main tools}

\subsection{Free factor systems, and automorphisms}

Let $G$ be a group. A finite collection of conjugacy classes of non-trivial subgroups  $\calG= \{ [H_1],  \dots , [H_r]\}$ will be called a free factor system of $G$ if there are representatives $H_1, \dots H_r$ of these conjugacy classes, and a subgroup $F_k$ of $G$, free of rank $k\geq 0$,  such that $G$ is a free product $G= H_1 * \dots *H_r*F_k$.
   The  \emph{Scott  complexity}  of this decomposition, and by extension of the free factor system $\calG$ of $G$,  is the pair $(k,r)$.  %

The set $\calT_\calG $ is the space of all metric $G$-trees whose vertex stabilizers are the conjugates of the groups $H_i$,  up to equivariant isometry. This space contains the (simplicial) Bass-Serre tree for the proposed decomposition. In a tree $T$ in   $\calT_\calG $, we say that a vertex $v$ is a free vertex if its stabilizer is trivial. If its stabilizer is a group whose conjugacy class is in $\calG$, we say that it is a non-free vertex.     See \cite{GL_out} for further references on this space of trees.

The subset $Hyp(\calG)$ of $G$ consists of all elements that are not conjugate into any of the subgroups $H_i$.  Those elements are said to be hyperbolic relative to $\calG$. They define loxodromic isometries of  the trees in $\calT_\calG$. We will say that elements conjugated to a subgroup $H_i$ are elliptic relative to $\calG$. They define elliptic isometries of the trees in  $\calT_\calG$.

On the level of automorphisms,  ${\rm Aut} (G, \calG)$ denotes the group of all automorphisms of $G$ that   preserve the conjugacy classes of each $H_i$.

Given any such automorphism $\phi$, and 
given any tree $T$ in $\calT_\calG$,  by \cite[Lem. 4.2]{SFAM} one can construct a  continuous map $f_\phi: T\to T$ that is $G$-equivariant with respect to $\phi$ in the following sense: for all $x\in T$ for all $g\in G$, $f(gx) = \phi(g) f(x)$. Such a map is called a topological realisation of $\phi$, and we say  that it represents $\phi$.  %

Let us now compare two differerent free factor systems.  If $G\simeq  H'_1 * \dots *H'_{r'}*F_{k'}$ for some other subgroups $H'_j$ and $F_{k'}$, and if  for each $i\in \{1, \dots ,r'\}$ there exists  $j\in \{1,\dots ,r\}$ and $g\in G$   satisfying $H'_i\leq gH_jg^{-1}$,  then one says that  the free factor system $ \calH'= \{[H'_1],...[H'_{r'}]\}$ is lower than $\calH= \{[H_1],...[H_r]\}$.  It is strictly lower if moreover $ \calH $ is not lower than  $\calH' $.     For instance, any Grushko's decomposition of $G$ provides a lowest free factor system for this order.   

If  $\calH'= \{[H'_1],\dots, [H'_{r'}]\}$ is strictly lower than $\calH=\{[H_1],\dots, [H_r]\}$, then  one may endow each $H_i$ with a free factor system inherited by its action on a tree of $\calT_{\calH'}$. Denote by $(k_i, r_i)$ the Scott complexity of this free factor system for $H_i$. The integer $r_i$ counts the number of conjugates of the $H'_j$ in its decomposition, and the integer $k_i$ counts the rank of the free group in its decomposition. Writting the identity $G= H_1*\dots H_r*F_k$ allows to show that $\sum_{i=1}^{r} r_i = r'$ and $( k +\sum_{i=1}^{r}  k_i) = k'$. This shows that, for the lexicographical order,  for all $i$, one has $(k_i, r_i) <(k', r')$ for all $i$, and that $(k, r)<(k', r') $. This is recorded in the next statement.

\begin{lem} \label{kurosh-scott} If the free factor system  $\calH'$ is strictly lower than $\calH$, then its Scott complexity is strictly larger than the complexity of $\calH$. If $[H] \in \calH$, then the Scott complexity of its free factor system induces by $\calH'$ is  strictly lower than the complexity of $\calH'$. 
\end{lem}

\subsection{Irreducibility, atoroidality, and twinned subgroups}

 We recall two equivalent  definitions of irreducibility of  automorphisms of a free product, proposed by Francaviglia and Martino \cite{SFAM}. They generalise the case of automorphisms of a free group.

\begin{defn}  %
\cite[Def. 8.1]{SFAM}
Let $ G$ be a group with a free factor system $\calG$. 
Let $\phi\in {\rm Aut}(G,\calG)$ and $T\in \calT_\calG$, and $f:T\rightarrow T$ representing $\phi$.

 We call $f$ \emph{irreducible}, if for every proper subgraph $W$ of the tree $T$ that is $G$-invariant and $f$-invariant, the quotient graph $G\backslash W$ is a collection of isolated subtrees with at most one non-free vertex. We say $f$ is \emph{fully irreducible} if for any integer $i>0$, $f^i$ is irreducible.

We say $\phi\in Aut(G,\calG)$ is $\calG$-fully irreducible if for any  $T\in \calT_\calG$, and for any  $f:T\rightarrow T$ representing $\phi$, the map $f$   is   fully irreducible. 
\end{defn}

This is equivalent to the following (see \cite[Def. 8.2, Lem. 8.3]{SFAM}):
\begin{defn}(Irreducible automorphisms relative to a free factor system)
Given a group $G$ with a free factor system $\calG$, 
  an automorphism $\phi\in Aut(G, \calG)$ is  irreducible (relative to $\calG$)     %
if $\calG$  is a maximal    %
proper free factor system that is invariant under $\phi$.  
\end{defn}

Let us make the following observation.  

\begin{lem} \label{lem-existence-fully-irre-ffs}
Let $G$ be a group with a free factor system $\calG$, and let $\phi \in {\rm Aut}(G, \calG)$. 
 
There exist a proper free factor system  $\calG'$ of $G$, preserved by some power of $\phi$, such that this power of $\phi$  is fully irreducible with respect  to $\calG'$.
\end{lem}

\begin{pf}
Denote $\calG$ by $\calG_1$, and $m_1=1$. Recursively, for all $n$, we construct $\calG_{n+1}$ from $\calG_{n}$, and $m_{n+1}$, a multiple of $m_n$,  as follows. 

 If $\phi^{m_n}$ is  fully irreducible relative to $\calG_n$, then $\calG_{n+1} = \calG_n$, and $m_{n+1}=m_n$. 
If $\phi^{m_n}$ is not  fully irreducible relative to $\calG_n$, some of its power $\phi^{m_{n+1}}$ preserves a strictly larger proper free factor system $\calG_{n+1}$, which hence has a strictly lower Scott complexity, by Lemma \ref{kurosh-scott}.  Since any lexicographically decreasing sequence in $\mathbb{N}\times \mathbb{N}$ is eventually constant, there is $n$ such that $\calG_n = \calG_{n+1}$.   This proves the lemma.

\end{pf}

Generalising a similar notion in free groups, we have the following. 

\begin{defn} (Atoroidal automorphisms)
We say $\phi\in Aut(G,\calG)$ is \emph{atoroidal}, if for any $g\in Hyp(\calG)$, and for any positive integer $n$, $[\phi^n(g)]\neq [g]$.
\end{defn}

We will need the following related notion.

\begin{defn}(Nielsen and pre-Nielsen paths)

Consider a group $G$ with a free factor system $\calG$, a tree $T\in \calT_\calG$, an automorphism $\phi\in {\rm Aut}(G,\calG)$, and a map $f: T\rightarrow T$ representing $\phi$. 

 A reduced path $\rho$ in $T$ is called a \emph{Nielsen} path if there exist an exponent $n\geq 1$ such that, for some $g\in G$,  the path $[f^{n} (\rho)]$ obtained by  $f^n(\rho)$ after path reduction, is equal to $g\rho$. A reduced path $\rho$ is called \emph{pre-Nielsen} if there exist an exponent $M>0$ such that $f^M(\rho)$ is Nielsen.
\end{defn}

Observe that even if $\phi$ is atoroidal, there can be Nielsen paths: they do not map on closed loops in $G\backslash T$.

We also introduce a related notion.

\begin{defn}\label{def;twinned}
Let $G$ be a group, $\calG$ be a free factor system, and $\phi \in {\rm Aut} (G, \calG)$. We  say that two different subgroups  $H, K$, such that $[H]$ and $[K]$ are preserved by $\phi$,  form a \emph{twinned pair of subgroups}  for $\phi$, if there exists $g\in G$ and an integer $m\geq 1$ such that   $\phi^m(H) = gH g^{-1}$ and $\phi^m(K) = gKg^{-1}$.
\end{defn}

  Let us underline that $[H]$ and $[K]$ are possibly (but not necessarily) equal. The groups can be in $\calG$, or in some other collection of preserved conjugacy classes of subgroups.

\subsection{Train Tracks Maps}

\subsubsection{Definitions and existence}

In \cite{BH} Bestvina  and Handel have defined a particular class of maps from a tree to itself, that is a cornerstone of the study of maps realising automorphisms.

In an oriented graph, let us denote by $i(e)$ the inital vertex of an oriented edge $e$.

\begin{defn}(Train track structure, legal turn, and legal paths) \cite{BH}

Given a graph $X$, an ordered pair $(e_1,e_2)$ of oriented edges such that $i(e_1) = i(e_2)$  is called a \emph{turn} (at the vertex $i(e_1)$).  A trivial turn is a turn of the form $(e,e)$.

A \emph{train track structure}  (or a gate structure) on a $G$-tree $T$ is a $G$-invariant equivalence relation on the set of oriented edges  at each vertex of $T$, with at least two equivalence classes  at each vertex.

Each equivalence class of oriented edges is referred to as a \emph{gate}.

In a gate structure, a turn is said to be \emph{legal} if the two oriented edges are in  different equivalent classes. A reduced path is said to be legal if all its turns  are legal.
\end{defn}

To describe a gate structure, it is enough to specify which turns are legal (or illegal).

An important example of gate structure is the one given as follows (and this is the one we will use). Consider $T$ and $T'$ two $G$-trees, and  a  map $f:T\to T'$ which is equivariant, and piecewise linear (linear, non constant, on edges). Define the gate structure on $T$ induced by $f$ as follows. Declare that a turn $(e_1, e_2)$ is illegal if $f(e_1)$ and $f(e_2)$ share their first edge in $T'$. It is easy to check that this defines an equivalence relation on the oriented edges issued from a same vertex, and that it is invariant for $G$, by equivariance of $f$.

In this construction, it is obvious that any legal turn is sent by $f$ on a pair of paths whose first edges define a non-trivial turn (by abuse of language we say that  any legal turn is send by $f$ on a  non-trivial turn). However, if $T'=T$, in principle, a legal turn could be sent on an illegal turn, and in that case $f^2$ would send a legal turn to a trivial turn. This is not a pleasant situation, and motivates the following.

\begin{defn}(Train track maps)
Given $T\in \calT_\calG, \phi\in Aut(G,\calG)$, and given $f:T\to T$  a piecewise linear $G$-equivariant map  (linear, non constant, on edges) realising $\phi$, we say that $f$ is a \emph{train track map} if, for the gate structure it defines,
\begin{itemize}
  \item $f$ maps edges to legal paths;
  \item if $f(v)$ is a vertex, then $f$ maps legal turns at $v$ to legal turns  at $f(v)$.
\end{itemize}

\end{defn}

One of the main results of Francaviglia and Martino in  \cite{SFAM}, is the following.

\begin{thm} \cite[Thm. 8.18]{SFAM}

If $\phi\in Aut(G,\calG)$ is irreducible, then there exist $T\in \calT_\calG$ and $f:T\rightarrow T$ representing $\phi \in Aut(G)$, such that $f$ is a train track map.
\end{thm}

Note also the following useful fact: by \cite[Lem. 8.20]{SFAM}, if $f:T\to T$ is a train track map representing $\phi$, then $f^k$ is a train track map representing $\phi^k$.

\subsubsection{Growth  rate of edges}

The metric point of view on train tracks is facilitated by the following. 

\begin{lem}(see also \cite[Lem. 8.16]{SFAM}) \label{lem;stretching_factor_Scott}

If $f:T\rightarrow T$ is a (piecewise linear) train track map representing a fully irreducible automorphism $\phi$, and if $T$ has at least two orbits of edges,  then there is a rescaling of each orbit of edge of $T$ such that, for this metric,  every edge is stretched by the same factor by $f$. More precisely, there is a constant $\lambda>1$ such that $l_T(f(e))=\lambda l_T(e)$ for all edge $e$ in $T$. 
\end{lem}

Such a constant $\lambda$ is called the growth rate of the train track map $f$. 

\begin{pf} \cite[Lem. 8.16]{SFAM} establishes the existence of  a rescaling to obtain the statement with $\lambda \geq 1$.  We show that under the extra assumption of the lemma, that $T$ has at least two orbits of edges, one has $\lambda >1$.

Assume  that this factor $\lambda$ is $1$. Then, partition the edges $\{e_1, \dots, e_s\}$ in subsets $E_1,\dots, E_{s'}$ of edges of equal length, from the shortest to the longest.   The set of the $G$-orbits of the edges in $E_{s'}$ (the longest) is permuted by $f$, otherwise one of them is not in the image of $f$. Hence, by iterating the argument, for all $i$ the set of $G$-orbits of the edges in $E_{i}$  is permuted by $f$. There is therefore $m$ such that $f^m$ sends any edge to one of its image by an element of $G$. 

 If, in $T$, there is an arc between two vertices fixed respectively  by $H_i$ and $H_j$, and whose all other vertices are free vertices (i.e. have trivial stabiliser), then the free factor $H_i*H_j$ is sent on of conjugate of itself by $\phi^m$.  Since $\phi$ is assumed to be fully irreducible, this means that either $G= H_1*H_2$ (and $T$ has only one orbit of edge), or $\calG$ has at most one free factor, $H_1$. We thus place ourselves in the later case: its Scott complexity is $(1,k)$.   

 Consider then  in $T$, a shortest  arc between two different points $v_1, v_2$ of the same orbit, and consisting only of free vertices, except possibly $v_1$ and $v_2$. Let us say that $v_2= gv_1$.  The  image in $G\backslash T$ of this arc has to be a simple loop (otherwise a strict subloop is suitable).  We may assume, after a possible translation, that the stabiliser of $v_2$ is either trivial or $H_1$. Thus, $\phi^m(g) = gx$ with $x$ in this stabilizer. It follows that $\phi^m$ preserves $\langle Stab(v_2), g\rangle$, which is either $\langle g \rangle$ or $H_1* \langle g \rangle$. In either case, full irreducibility of $\phi$ forces $G$ to be $H_1*\langle g\rangle$. Hence,   $\calG$ has a Scott complexity $(1,1)$,  and $T$ has one orbit of edges.  

 We have proved that if $T$ has at least two orbit of edges, then $\lambda >1$.

\end{pf}

\subsection{Coning off subtrees}

In this subsection, we consider $G$ a group, with a free factor system $\calG$, %
and a $G$-tree $T$ in $\calT_\calG$. 

Let us assume that we are given a collection of conjugacy classes of subgroups of $G$,  $\mathcal{P}= \{[P_1], [P_2], \dots, [P_q]\}$, invariant under $\phi$ (in the stronger sense that each $[P_i]$ is invariant). For each $i$, let $T_i$ be the minimal subtree of $T$ 
 invariant under $P_i$. 

In this situation, one considers the cone-off tree $\dot T$ of $T$,  adding one vertex $v_{ gP_ig^{-1}}$ for each coset $gP_i$ and adding edges from $v_{ gP_ig^{-1}}$ to all free vertices of $gT_i$ (by subdividing edges in an orbit of $T_i$, if necessary, we may always assume that there are such free vertices).       
We choose the length of these new edges to be smaller than half the minimal length of edges in $T$, and smaller than $1/2$. The length of a segment $\sigma$ will be denoted $l_{\dot T} (\sigma)$.

We say that $\calP$ is above $\calG$ if any group whose conjugacy class is in $\calG$ is contained in a group whose conjugacy class is in $\calP$. Recall that $\calP$ is malnormal if for each indices $i,j$, and each $g\in G$, if $P_i\cap gP_jg^{-1}$ is non-trivial, then $i=j$ and $g\in P_i$. 

\begin{pro} %
  Assume that $\calP$ is a finite collection of conjugacy classes of subgroups of $G$ that is above $\calG$. Assume that for all $i$, the action of   $P_i$ on $T_i$ is cofinite, and that the collection $\calP$ is malnormal. Then $G$ is hyperbolic relative to $\calP$, and $\dot T$ is hyperbolic. 
\end{pro}

\begin{pf} The assumption on cofiniteness of $P_i$ on $T_i$ ensures that $P_i$ is relatively quasi-convex in $G$ relative to $\cal G$. One can then apply Yang's \cite[Thm 1.1]{Wenyuan_Crelle} to get the conclusion.

  \end{pf}

  We will say that $\calP$ is \emph{hyperbolically coning-off} $\calG$ if it satisfies the four assumptions of the Proposition: $\calP$ is a finite collection of conjugacy classes of subgroups of $G$, $\calP$ is malnormal, $\calP$ is above $\calG$, and for all $i$, the action of   $P_i$ on $T_i$ is cofinite.

  We then define a reduction in $\dot T$ of a path of $T$ as follows. Call a path in $T$ collapsible if its end points are neighbors of a same vertex $v_{gP_i}$, and its collapse is the two-edge path between its end points that goes through $v_{gP_i}$. For an arbitrary reduced path in $T$, choose a maximal collection of maximal collapsible subsegments, and replace each of them by its collapse.
  The reduction of a reduced path in $T$ is a uniform quasigeodesic in $\dot T$ (\cite[Prop. 2.11]{DMj}). 

  If $f: T\to T$ is a map representing an automorphism $\phi$ that preserves every conjugacy class of groups in $\calP$, then $f$ induces $\dot f :\dot T\to \dot T$ by sending, for each $P$ (with $[P]$ in $\calP$),  the vertex fixed by $P$ on the vertex fixed by $\phi(P)$ (and edges according to the action). Observe that if $f$ is  a train track map, then every turn at a vertex fixed by $P$ is a legal turn for $\dot f$, and is sent on a legal turn.

\subsection{Polynomially growing subgroups}\label{sec;321} %

Consider $F$  a free group, and $\phi$ and automorphism of $F$. We briefly discuss some material covered in the preliminary section of \cite{Lev09}, to which the reader is warmly refered.

A subgroup $F_0$ of $F$ is said \emph{polynomially growing}  for  $\phi$ if for every $g\in F_0$, the length $|[\phi^n(g)]|$ of a cyclically reduced element in the conjugacy class $[ \phi^n(g) ]$ is bounded above by a polynomial in  $n$. 
If $F$ itself is polynomially growing, one says that $\phi$ is a polynomially growing automophism.

Recall that the group $Aut(F)$ contains as a normal subgroup the group of all inner automorphisms $Inn(F)$, which are conjugations by elements of $F$, and the quotient  $Aut(F)/Inn(F)$ is the outer automorphism group of $F$.

For an  outer automorphism $\Phi\in Out(F) $, we say that a
subgroup $F_0$ of $F$ is polynomially growing  for $\Phi$ if there is $\phi$ and automorphism in the class of $\Phi$,  for  which $\phi(F_0) = F_0$ and $\phi|_{F_0}$  is  polynomially growing.

In \cite[Prop 1.4]{Lev09} Levitt proves that for any outer
automorphism $\Phi$ of a free group $F$, there is a finite family of
finitely generated  subgroups of $F$, that are polynomially growing for $\Phi$,
such that all polynomially growing 
subgroups for $\Phi$ are conjugated into one of them. These are maximal polynomially growing subgroups for $\Phi$.  Levitt also proves in the same reference, that the family of their conjugates is malnormal in sense recalled in Section \ref{sec;Mappingtori}.

We now adapt this setting to the case of free products. 

Let $G$ be a group and $\calG$ be a free factor system for $G$. Let $\phi\in {\rm Aut} (G,\calG)$ be a fully irreducible automorphism. Let $f:T\to T$ be a train track map representing $\phi$ on a tree $T\in \calT_\calG$.

We define polynomial growth for  $\phi$ in $T$ as follows.   %
 An element $h$ of $G$ is said to be \emph{polynomially growing} (or to have polynomial growth)  for $\phi$ in $T$,  if the translation length of $\phi^n(g)$ in $T$ is bounded above by a polynomial in $n$.    A subgroup $P$ of $G$ is said polynomially growing  for   $\phi$ if there exists $g\in G$, and $k\geq 1$  such that for $\psi = {\rm ad}_g\circ \phi^k$ one has $\psi( P) = P$ and such that  for every $h\in P$, the translation length of $\psi^n(g)$  is bounded above by a polynomial in  $n$.  Note that this property is invariant if one considers a power of $\phi$ instead of $\phi$. 

Observe that if the stretching factor $\lambda$ is $1$,  then the whole group $G$ is polynomially growing for $\phi$. We assume now that the stretching factor $\lambda$ is $>1$, or in view of Lemma \ref{lem;stretching_factor_Scott}, that $T$ has more than one orbit of edges under the $G$-action, or in other words, that the Scott complexity is different from $(0,2)$ or $(1,1)$.

We now recall the construction of a limiting $\mathbb{R}$-tree, found in \cite[\S 2.5]{Gab_et_al}, to which we refer for the details  (our setting is slightly eased by the assumption that $\phi$ is fully irreducible, that is $\tau' = \emptyset$ in the notation of  \cite[\S 2.5]{Gab_et_al}).   One may endow $T$ with a sequence of actions of $(G,\calG)$, by precomposing by powers of $\phi$, and rescaling by the factor $\lambda$.  %

One may then go to an ultralimit of these metrics and obtain an $\mathbb{R}$-tree $T^\infty$ endowed with an action of $(G,\calG)$ which is non-trivial and minimal (see \cite[Lem. 2.7]{Gab_et_al}),
and with an homothety of scale $\lambda$.
This action has been considered in several places, we indicate to the reader \cite[\S 2.5]{Gab_et_al} as a reference.  We gather a few useful facts about it.

By \cite[Lemma 2.8]{Gab_et_al} (its proof applies without change in our case), this action has trivial arc stabilizers. 
Let us observe that any polynomially growing subgroup of $G$ fixes a point in $T^\infty$. Indeed in the rescaling process, all its elements have translation length going to $0$, and, even if we do not know yet that these subgroups are finitely generated, we may still deduce that they all fix a unique point, since arc stabilizers are trivial. By  \cite[Prop. 4.4]{Horb_Boun},  point stabilizers in $T^\infty$ have strictly smaller Kurosh rank than $G$. 

By \cite[Coro. 4.5]{Horb_Boun}, there are finitely many orbits of branch points in $T^\infty$. It follows that  if $G_0$ is a point stabilizer in $T^\infty$, there is $k\geq 1$, and $g_0\in G$ such that ${\rm ad}_{g_0} \circ \phi^k$  preserves $G_0$.   If moreover $P <G_0$ is a polynomially growing subgroup of $G$ for $\phi$ on $T$, then on the minimal subtree $T_0$ of $G_0$ in $T$,  the subgroup $P$ is polynomially growing for ${\rm ad}_{g_0} \circ \phi^k$.

Following \cite{Lev09}, we discuss the next five properties.

\begin{pro}\label{prop;Levitt_four_points}
  If $P_0$ is a polynomially growing subgroup of $G$ for $\phi$ on $T$, (resp. if $h$ is a polynomially growing element for $\phi$ in $T$), then there exists a unique maximal polynomially growing subgroup $P$ of $G$ for $\phi$ on $T$, such that $P_0<P$ (resp. such that $h\in P$).

  The collection of maximal polynomially growing subgroup of $G$ for $\phi$ on $T$ is malnormal.

 Any maximal polynomially growing subgroup of $G$ for $\phi$ on $T$ acts cofinitely on its minimal subtree  in $T$ and has Kurosh rank smaller than the rank of $G$ (strictly if it is not $G$). 
  
  If  $P$ and $P'$ are two maximal polynomially growing subgroup of $G$ for $\phi$ on $T$, such that there exists $k\geq 1$ for which $\phi^k(P) = P, \phi^k(P') = P'$, then $P=P'$.
  
  There are only finitely many conjugacy classes of maximal polynomially growing subgroups for $T$, in $G$.
\end{pro}

\begin{pf}
 If the rank is $1$ (or more generally if the Scott complexity is $(2,0)$ or $(1,1)$),  then $G$ itself is polynomially growing, and the statements are obvious. %

The five statements are obtained by induction on the Kurosh rank of $G$. We may assume that the statement is established for groups of Kurosh rank smaller than $G$, and that $G$ is not polynomially growing itself. 
Recall that by  \cite[Prop. 4.4]{Horb_Boun},  point stabilizers in $T^\infty$ have strictly smaller rank than $G$ (strictly, since they are different from $G$). 

First $P_0$ (resp. $h$) is elliptic in $T^\infty$,   hence is contained in a  point stabilizer in $T^\infty$, unique (since arc stabilizers are trivial),  call it $G_0$, as is any polynomially growing subgroup containing $P_0$ (resp. containing $h$). The induction hypothesis thus proves the first point.  

If $P, P'$ are two maximal polynomially growing subgroups that intersect, they must be contained in a same point stabilizer in $T^\infty$, and the induction hypothesis shows that $P=P'$. For malnormality it remains to check that any such $P$ is its own normalizer. If $g$ normalises $P$, then $g$ must be in $G_0$ the point stabilizer on $T^\infty$ containing $P$, and by induction, $g\in P$. The second point is proved. 

The rank statement in the third point also follows by induction, since by induction $P$ has Kurosh rank smaller than that of $G_0$ containing it, and $G_0$ has Kurosh rank strictly smaller than $G$. The statement on cofiniteness on its minimal tree follows similarily. First, $G_0$ acts on its minimal subtree, and inherits the structure of graph of groups with trivial edge groups from the quotient graph. Since its rank is finite, the quotient graph is finite, and $G_0$ acts cofinitely on its minimal subtree. By induction then, $P $ acts cofinitely on its  minimal subtree inside that of $G_0$, which is its minimal subtree in $T$.

For the fourth point, assume that  $\phi^k(P) = P, \phi^k(P') = P'$, then $P$ and $P'$ must fix the same point in $T^\infty$, since $\phi$ produces a dilatation of factor $\lambda$ in $T^\infty$. Then the statement follows by induction.

Finally for the fifth point, first by \cite[Coro. III.3]{GabL}, there are only finitely many $G$-conjugacy classes of stabilizers of points in $T^\infty$. In each of them, we may use the induction hypothesis to find that there are finitely many conjugacy classes of maximal polynomially growing subgroups.

\end{pf}

Let now $\calP$ be  the collection of conjugacy classes of maximal polynomially growing subgroups for $\phi$ on $T$, and $\dot T$ the coned-off tree of $T$ for this collection, as in the previous subsection.
  
We say that $\phi$ is atoroidal in $\dot T$ if for all $g\in G$ hyperbolic on $\dot T$, and all $m>0$, $\phi^m(g)$ and $g$ are not conjugate in $G$.  %

\begin{pro} \label{prop;atoroidal_P}  If $\calP$   %
  is the collection of conjugacy classes of maximal polynomially growing subgroups for $\phi$ on $T$, %
  then $\phi$ is atoroidal in $\dot T$.  %
\end{pro}  

\begin{pf} If $\phi$ was not atoroidal with respect to $\calP$, there is an element of polynomial growth (of degree 0 !) that is not in a group of $\calP$, hence $\calP$ does not cover the polynomial growth.

  \end{pf}

Recall  that by assumption, $\phi$ is fully irreducible for $\calG$, and that $f:T\to T$ is a train track map.

\begin{pro} (Transversality of legal paths) \label{prop;transversality}

Let $G$ be a group, $\calG$ a free factor system, $\phi\in {\rm Aut} (G, \calG) $ be fully irreducible. Let $f: T\to T$ be a train track map representing $\phi$, for some $T\in \calT_\calG$. 

Let $\dot T$ be the cone off of $T$ over the family of maximal polynomially growing subgroups of $G$ in $T$.  %
  
There exists $A>1$ such that 
if $\rho$ is a legal path in $T$, then the length of its reduction in $\dot T$ is larger than $\frac{1}{A} l_T(\rho)$.   
\end{pro}

\begin{pf}
Assume the contrary, that for all $n$ there exists $\rho_n$, a legal path, for which the length of its reduction in $\dot T$ is less than $\frac{1}{n} l_T(\rho)$.  Then  $\rho_n$ must contain a long subsegment $\tau_n$ contained in a translate of some $T_i$. Since the action of $H_i$ on $T_i$ is co-finite, $\tau_n$ must pass several times through free vertices in the same $P_i$-orbit. This means that $\tau_n$ contains $\eta_{n,1} \subset  \dots \subset  \eta_{n,3}$ subpaths with same initial point (that is a free vertex) that are legal, and that are the fundamental segments of  elements $\gamma_{n,j}$, conjugated in some $P_i$, that are loxodromic in $T$.   Denote by $e_{n,j}$ the last edge of $\eta_{n,j}$. Now, consider the turn made by $\gamma_{n,j}^{-1} e_{n,j}$, and   the (common) first edge $e_0$ of $\eta_{n,j}$. If one of these turn is legal, then the whole axis of $\eta_{n,j}$ is a legal path, and the growth of the conjugacy class of $\phi^m(\eta_{n,j})$ is exponential, contradicting its belonging to $P_i$.

If none of them is legal, then consider the path $\mu$ between the end point of $\eta_{n,1}$ and the last point of $\eta_{n,2}$. It is legal, and since $\gamma_{n,1}^{-1} e_{n,1}$ and $\gamma_{n,2}^{-1} e_{n,2}$ both are in the same gate than $e_0$, the turn made by $\gamma_{n,1} \gamma_{n,2}^{-1} e_{n,2}$ and the first edge $e_1$ of $\mu$ must be legal, as is the turn of $e_{n,1}$ and $e_1$ by assumption. We can use the same argument than before to get the desired contradiction. 
\end{pf}

\subsection{Hyperbolicity of an automorphism}

Let $G$ be a finitely generated group with a generating set $S$ and Cayley graph $\Gamma_S(G)$. Let $\Lambda$ be a set and let $\mathbb{H}=\{H_i\}_{i\in \Lambda}$ be a family of subgroups $H_i$ of $G$;

The $\mathbb{H}$-\emph{word metric} $|\cdot|_{\mathbb{H}}$ is the word-metric for $G$ equipped with generating set $S_{\mathbb{H}}=S\cup({\cup_{i\in \Lambda} H_i})$.

The following definition, from \cite{GW}, might appear technical. It extends a similar definition of hyperbolicity of automorphisms.  We will only use it for free products endowed with a free factor system, so the reader is free to restrict the definition to this case.

\begin{defn}\label{def-rel-hyp-aut}(Relatively hyperbolic automorphisms)

Let $G$ be a group with a generating set $S$. Let $\Lambda$ be a set and let $\mathbb{H}=\{H_i, i\in \Lambda\}$ be a family of subgroups $H_i$ of $G$ such that each $H_i$ is its own normalizer.  Let $\calH= \{  [H_i], i\in \Lambda \}$.  An automorphism $\phi\in Aut(G,\mathcal{H})$ is \emph{hyperbolic relative to} $\mathbb{H}$ (or in short, \emph{relatively hyperbolic}) if it satisfies the following:

there exist $\lambda>1$, $M,N \geq 1$, such that for any $g\in G$ with $|g|_\mathbb{H}\geq M$, the inequality holds:
$$\lambda|g|_\mathbb{H}\leq \max\{|\phi^N(g)|_\mathbb{H},|\phi^{-N}(g)|_\mathbb{H} \}$$
\end{defn}

Recall that for a path $\rho$ in a metric tree $T$, $l_T(\rho)$ is its length. Assume now that $G$ is a free product.
For a hyperbolic element $g_0$ in $G$, a \emph{fundamental segment} in $T \in \calT_\calG$ is a segment in its translation axis that starts and ends at free vertices $x, y$ (i.e. with trivial stabilizers),   such that $g_0x=y$.  Up to subdividing edges of $T$, such a fundamental segment for $\gamma$  always exists.

\begin{defn} \label{def;hyp_pair}
Let $G$ be a finitely generated group, and $\calG$ be a free factor system, and $\calP$ a collection hyperbolically coning off $\calG$.

Let $T, T'$ be metric trees in $\calT_\calG$, and $\dot T, \dot T'$ their cone-off over $\calP$. 

 Let $\alpha: \dot T\rightarrow \dot T'$ and $\alpha': \dot T'\rightarrow \dot T$ be $G$-equivariant Lipschitz maps between these spaces.

 Let  $\phi\in Aut(G,\calG)$, preserving $\calP$, and consider   $\dot f:\dot T\rightarrow \dot T$, $\dot f': \dot T'\rightarrow \dot T'$ induced on the cone-off by  maps representing $\phi\in Aut(G,\calG)$ and $\phi^{-1}\in Aut(G,\calG)$ respectively.

If there exist a natural number $M>0$ and a real number $\lambda>1$ 
such that  whenever $\sigma$ is either a path between two different non-free vertices of $\dot T$, or a fundamental segment in $\dot T$ for a hyperbolic element of $G$,   and  
and  $\sigma'$ the  pull-tight of its image by $\alpha$     
one has:

$$\lambda l_{\dot T}(\sigma)\leq \max\{l_{\dot T}[\dot f^{M}(\sigma)], l_{\dot T}([\dot f'^{M}(\sigma')])\}, $$

then the pair $(\dot f, \dot f')$ is refered to as a hyperbolic pair for $G,\calG, \phi$.

\end{defn}

\begin{rmk}\label{rem;pass}
Notice that, given $G$ and $\calH=\{[H_i], i=1,\dots, p\}$ a free factor system,    if $\mathbb{H}=\{H_1,...,H_p\}$, and if  $T$ is in $\calT_\calG$, then  the $\mathbb{H}$-word metric is quasi-isometric to the distance induced by an orbit in $T$. Therefore, if $(f,f')$ is a hyperbolic pair, then $\phi$ it is relatively hyperbolic.
\end{rmk}

The relevance of hyperbolic automorphisms, for us,  is through the following combination theorem, which is proven by Gautero and Weidmann. For a definition of relatively hyperbolic group, we refer to \cite{Bow}.

\begin{thm} \cite[Coro 7.3]{GW} \label{GW-comination-theorem-full}
  
Let $G$ be a finitely generated group that is hyperbolic relative to a finite family of conjugacy classes of  infinite subgroups $\mathcal{H}=\{[H_1], \dots, [H_k]\}$. Let us write $\mathbb{H}$ for $\{[H_1], \dots, [H_k]\}$.

 Assume that   $\alpha\in {\rm Aut} (G, \mathcal{H})$ is hyperbolic relative to $\mathbb{H}$.  Assume that for all $i=1, \dots, k$, the element $g_i\in G$ is such that $g_i^{-1} \alpha(H_i) g_i = H_i$. 

Then the semi-direct product $G\rtimes_{\alpha}\mathbb{Z}$ is relatively hyperbolic  with respect to the subgroups $\langle H_i, t_\alpha g_i \rangle \simeq H_i \rtimes_{ {\rm ad}_{g_i} \circ \alpha}   \mathbb{Z}$, where $t_\alpha$ denotes the generator of the semidirect product associated to the automorphism $\alpha$.

\end{thm}

For brevity one calls the collection of conjugacy classes of  subgroups  $\{   [\langle H_i, t_\alpha g_i \rangle] \}$, the mapping torus of $\mathbb{H}_0$.

As an alternative of the previous theorem, one may also use  the following form of \cite[Thm 4.6]{MjR}, for the case of an HNN extension. We will explain how to use it instead, when we will need it.

\begin{thm} (Mj-Reeves, Particular case of \cite[Thm 4.6]{MjR}) \label{thm;MjR}

Let $\Gamma$  be the fundamental group of  a finite graph  of  relatively hyperbolic groups satisfying
(1) the qi-embedded condition, 
(2) the strictly type-preserving condition,
(3) the qi-preserving electrocution condition,
(4) the induced tree of coned-off spaces satisfies the hallways flare condition,
(5) the cone-bounded hallways strictly flare condition.
Then $\Gamma $ is  hyperbolic relative to the family $\calC$ of maximal parabolic subgroups.   
\end{thm}

\subsection{More mapping tori: global preservation of a collection of subgroups}  \label{sec;Mappingtori}

Consider the case where $G$ is a group, $\phi$ is an automorphism, and $\calG=\{[H_1], \dots, [H_p]\}$ is a malnormal collection of conjugacy classes of subgroups of $G$, globally preserved by $\phi$. Up to taking a power we retrieve the preferred case in which each $[H_i]$ is preserved by $\phi^{p!}$. In this subsection, we argue that, in the interesting cases,  the relative hyperbolicity result that one may have for the mapping torus by $\phi^{p!}$ descend to the mapping torus by $\phi$.

 We first define mapping tori of groups in the collection $\calG$.

 If  $\phi$  preserves the conjugacy class of each $H_i$,  let us choose $g_i$ such that $\phi(H_i) =g_i H_i g_i^{-1}$. Note that by malnormality, $g_i$ is unique up to a multiplication on the right by an element of $H_i$.  Consider the semidirect product $G\rtimes_\phi \mathbb{Z}$, and call $t$ the generator of the factor $\mathbb{Z}$ that acts on $G$ as $\phi$. We call the mapping torus of the collection $\calG$  by $\phi$, the collection of conjugacy classes of subgroups $\langle H_i, tg_i\rangle$ in $   G\rtimes_\phi \mathbb{Z}$, for $i=1, \dots, p$. Note that $tg_i$ normalises $H_i$, and that these subgroups are uniquely defined (a change in the choice of $g_i$ does not change the subgroup).

If now $\phi$ preserves $\calG$ globally,  then there is $m$ (dividing $p!$) such that $\phi^m \in  {\rm Aut}(G,\calG)$.

In that case, for each $i$, there is $m_i$ for which $\phi^{m_i}$ preserves the conjugacy class of $H_i$. One may define the mapping torus of $H_i$ by $\phi^{m_i}$ as before (choosing an element $\gamma_i$ so that $t^{m_i}\gamma_i$ normalises $H_i$), and we declare that the mapping torus of the collection $\calG$ by $\phi$ is the collection of conjugacy classes of mapping tori of $H_i$ by $\phi^{m_i}$.

\begin{lem} \label{lem;maln_finite_distance} Let $G$ be a free group. If $H$ is a non-trivial, malnormal finitely generated subgroup in $G$ and if  $K$ is the image of $H$ by an automorphism of $G$, that is at bounded distance from $H$ in the word distance, then $K=H$. 
\end{lem}
\begin{pf}
Since $H$ is malnormal it is not a proper finite index subgroup of another subgroup of $G$. Since $K$ is an automorphic image, the same is true for $K$.  However, being finitely generated, they both are of finite index in the stabilizer of their common limit set in the boundary of the free group $G$. Therefore they are equal. 
\end{pf}

\begin{pro}\label{prop;phimtophi}
  With the notations above, assume that $G$ is free, and that the groups in $\calG$ are infinite.

  If $G\rtimes_{\phi^m}\mathbb{Z}$ is relatively  hyperbolic relative the mapping torus of $\calG$ by $\phi^m$, then $G\rtimes_{\phi}\mathbb{Z}$ is relatively hyperbolic relative to  the mapping torus of the collection $\calG$ by $\phi$. 
  \end{pro}

  \begin{pf}
    The group $G\rtimes_{\phi^m}\mathbb{Z}$ is a finite index subgroup of  $G\rtimes_{\phi}\mathbb{Z}$, therefore, by Dru\c{t}u's theorem on invariance of relative hyperbolicity by quasi-isometry \cite[Thm. 5.1]{Dru},  $G\rtimes_{\phi}\mathbb{Z}$ is relatively hyperbolic with respect to a collection of subgroups such that each is at bounded distance from a peripheral subgroup in $G\rtimes_{\phi^m}\mathbb{Z}$. Consider $Q$ a peripheral subgroup of $G\rtimes_{\phi}\mathbb{Z}$, and (possibly after conjugation) let  $\langle H_i, tg_i\rangle$ be the peripheral subgroup of $G\rtimes_{\phi^m}\mathbb{Z}$ that remains at bounded distance. First, $Q\cap G$ must be at bounded distance from  $H_i$, therefore equal, by the previous lemma.        
Second, if $z\in Z$ is not in $G$, it conjugates $H_i$ into some subgroup of $G$, therefore it must normalise $H_i$. It follows it must be of the form $(t^{m_i} \gamma_i)^s h$ for some $h\in H_i$. This ensures $Q$ is contained in the mapping torus of $H_i$ by $\phi^{m_i}$. Conversely, since $t^{m_i} \gamma_i$ normalises $H_i$ which is an infinite parabolic subgroup, it must be in the associated peripheral subgroup. This proves that $Q$ is the  mapping torus of $H_i$ by $\phi^{m_i}$, and it proves our proposition. 
    \end{pf}

\section{Relative Hyperbolicity of semi-direct products}
\subsection{Statements in full irreducibility}

We first state the following result,  an analogue of Brinkman's first result in \cite{PB}, and of a result of Bestvina-Feighn-Handel for free groups in \cite[Thm. 5.1]{BFH}.

\begin{thm}\label{thm;rel-hyp-fully-irreducible}
Let $G$ be a finitely generated group with a free factor system $\calG$, of Scott complexity $(k,p)$, different from $(1,1)$ and $(0,2)$.    Let $\phi\in {\rm  Aut}(G,\calG)$ be  fully irreducible and atoroidal. Assume that it has no twinned subgroups in $\calG$. 

 Then the semi-direct product $G\rtimes_{\phi}\mathbb{Z}$ is relatively hyperbolic, with respect to the mapping torus of $\calG$.     
\end{thm}   

We will actually prove the following variant, which immediately implies Theorem \ref{thm;rel-hyp-fully-irreducible},  by choosing $\calP$ to be equal to $\calG$.

\begin{thm}\label{thm;rel-hyp-fully-irreducible_dot}
  Let $G$ be a finitely generated group with a free factor system $\calG$, of Scott complexity $(k,p)$, different from $(1,1)$ and $(0,2)$.     Let $\phi\in {\rm  Aut}(G,\calG)$. %

  Let $\calP$ be  %
  a  $\phi$-invariant (in the strong sense that $\phi \in {\rm Aut}(G,\calP)$)  collection that is hyperbolically coning-off $\calG$.
  
  Assume that $\phi$ is fully irreducible with respect to $\calG$, atoroidal with respect to $\calP$, and has no twinned subgroups in $\calP$.

  Assume that  for a tree  $T\in \calT_\calG$ and a train track map $f:T\to T$ realizing $\phi$, there exists $A$ such that  if $\rho$ is a legal path in $T$ of length $l_T(\rho)$ in $T$, its reduction in $\dot T$ is of length at least $ \frac{l_T(\rho)}{A} $, and similarily for $\phi^{-1}$.

 Then, for some power $\phi^m$ of $\phi$, (with $m =1$ if $\calP=\calG$)  the semi-direct product $G\rtimes_{\phi^m}\mathbb{Z}$ is relatively hyperbolic, with respect to the mapping torus of $\calP$.     
\end{thm}

The last condition, about the existence of $A$  will be refered to as $\calP$ being transversal to legal paths for $\phi$.

Observe that if $\calP$ is the collection of maximal polynomially growing subgroups of $G$ in $T$, then the assumptions of the theorem are satisfied, by Propositions \ref{prop;Levitt_four_points}, \ref{prop;transversality} and \ref{prop;atoroidal_P}. We thus obtain from Theorem \ref{thm;rel-hyp-fully-irreducible_dot} the following corollary.

\begin{cor} \label{coro;23} Let $G$ be a finitely generated group with a free factor system $\calG$, of Scott complexity $(k,p)$, different from $(1,1)$ and $(0,2)$.     Let $\phi\in {\rm  Aut}(G,\calG)$  be  fully irreducible.

  Let $\calP$ be  the collection of conjugacy classes of maximal polynomially growing subgroups for $\phi$ on a (any) tree in $\calT_\calG$.  Then, there is a power of $\phi$ for which $\phi^m\in {\rm Aut} (G, \calP)$, and for which  the semi-direct product $G\rtimes_{\phi^m}\mathbb{Z}$ is relatively hyperbolic, with respect to the mapping torus of $\calP$.  
\end{cor}

 The proof of Theorem \ref{thm;rel-hyp-fully-irreducible_dot}  will take the next subsections, until the end of Section \ref{sec;proof_fully_irred}, where it will be formally given. %
 As discussed in the introduction, it will follow closely Brinkmann's and Bestvina Feighn and Handel's proofs \cite{PB,BFH}, with the additionnal difficulty of dealing with the lack of local finiteness of the trees  involved, and their cone-off.

\subsection{Growth in train tracks}\label{sec;21}

In this section, $G$ is a group, endowed with a free factor system $\calG= \{[H_1], \dots, [H_p]\}$, of Scott complexity different from $(2,0)$ or $(1,1)$, and a collection $\calP = \{[P_1], \dots, [P_q]\}$ hyperbolically coning-off $\calG$. We will denote $ H_{p+j} =P_j$ for convenience and unification of notations.

Recall that, if $\phi\in {\rm Aut} (G,\calG)$ is a fully irreducible  automorphism, that preserves $\calP$,    we have agreed (in Lemma \ref{lem;stretching_factor_Scott}) on a metric on a train track tree so that each edge is stretched by the train track map $f$ by a factor $\lambda >1$,  as the Scott complexity is different from $(2,0)$ or $(1,1)$.

\subsubsection{Angles}

The following introduces our tool for coping with the non-local finiteness of the trees in $\calT_\calG$ and their cone-off, an issue that already showed in the previous lemma.

\begin{defn}
For each $i \leq p+q$, fix a word metric $|\cdot|_i$ on $H_i$. Let $T\in \calT_\calG$ and  let $v_{H_i}$ be a vertex of $T$  or $\dot T$, that is fixed by $H_i$. Choose  $\calE_{v_{H_i}} = \{\varepsilon_0, \dots, \varepsilon_r\}$  a transversal of adjacent edges of $v_{H_i}$ for the action of its stabiliser. 

 For each pair of edges $e,e'$ adjacent to $v_{H_i}$,  the angle $Ang_{v_{H_i}}(e,e')$ is the word length of $g^{-1}g'$, where $g,g'$ satisfy that bot  $ge_1, g'e_2$ are in $ \calE_{v_{H_i}}$.

 For each $i$, for each $v\in Gv_{H_i}$, and for each pair of edges $e_1, e_2$ adjacent to $v$, define the angle $Ang_v(e_1,e_2)$ to be $Ang_{v_{H_i}}(g''e_1,g''e_2)$ where $g''$ is the element in $G$ such that $g''v=v_{H_i}$.

We say that a path is $\Theta$-straight if  angles between its consecutive edges are at most   $\Theta$.

\end{defn}

\begin{rmk}
Angles are well-defined. Indeed, the choice of $g$ and $g'$ is unique, because stabilizer of each edge (and thus of each edge in $\calE_{v_{H_i}}$)  is trivial.
In general, the choice of $g'' $ in the definition is not unique, but only differs from an element in $H_i$, and by the definition of angles at the vertex whose stabilizer is $H_i$, element in $H_i$ preserves the angle. 

We also notice that the angles safisfy a local finiteness: for a given edge $e_1$ with starting vertex $v$ and a given number $C>0$, there are only finitely many possible $e_2$ satisfying $Ang_v(e_1,e_2)<C$. This is easy to see as there are only finitely many edges (up to $G$-orbit) adjacent to $v$ and that there are only finitely many elements in $H_i$ whose word length is bounded by $C$.

Finally, we notice that  angles are $G$-invariant. 
\end{rmk}

\begin{lem}\label{lem-finitely-many-theta-equivalent-paths}
Let $\phi\in Aut(G,\calG)$ be an automorphism of $G$,    $f:T\rightarrow T$ be a map representing $\phi$, and $\dot f: \dot T\to \dot T$ the induced map on $\dot T$. Then for any $\Theta_1>0$, there exist $\Theta_2>0$, such that for any pair of edges $e_1, e_2$ starting from a vertex $v$ with $Ang_v(e_1,e_2)>\Theta_2$, we have that $Ang_{f(v)}(\dot f(e_1,  \dot f(e_2)))>\Theta_1$.
\end{lem}

\begin{pf} Assume the contrary, let $v_n, e_n, e'_n$ such that $Ang_{v_n}(e_n,e'_n)>n$ but for which  $Ang_{f(v_n)}(\dot f(e_n, \dot f(e'_n))) \leq \Theta_1$. After translation, and extraction, we may assume that $v_n$ is constant (we'll denote it by $v$), and stabilised by $H_i$ for some $i\leq p+q$, and that $e_n$ is also constant (denoted by $e$).  Let $h$ and $j$ such that $\dot f(v)$ is stabilised by $h^{-1}H_jh$. Since $f$ represents $\phi$, the automorphism ${\rm ad}_{h} \circ \phi$ induces, by restriction to $H_i$ an isomorphism to  $H_j$, hence a quasi-isometry for their word metrics.   However, that $Ang_{v}(e,e'_n)>n$  indicates that there is $h_n$,   a sequence of elements of $H_i$ going to infinity, such that $h_n^{-1}e'_n$  remains among finitely many edges. After extraction it is constant, $e'$. The image of $e'_n$ by $\dot f$ correspond to the images of $\dot f(e')$ by $\phi (h_n)$ hence is going to infinity in angle from $\dot f(e)$. This is a contradiction. 
\end{pf}

From the above lemma, we have:

\begin{lem}\label{e-g-fs-L}
Let $G$ be a group with a free factor system $\calG$, a collection $\calP$ that is hyperbolically coning-off $\calG$ as before, an automorphism $\phi\in {\rm Aut}(G,\calG)$ preserving $\calP$. 
Assume that the Scott complexity of $(G, \calG)$ is different from $(1,1), (0,2)$,  and that $\phi$ is    
atoroidal for $\calP$. Let $T\in \calT_\calG$, and  $f:T\to T$  
 representing $\phi$ on $T$. Let $\dot f: \dot T\to \dot T$ the induced map on the cone-off.

 Then for any given $h\in G$ hyperbolic on $\dot T$, any fundamental segment $\tau$ of $h$ in $\dot T$, and for any $C>0$, there is an integer $N>0$ such that $l_{\dot T}(\dot f^N(\tau))>C$.

Moreover, if $\tau$ is a path in $\dot T$ between two non-free vertices, and if $\phi$ has no twinned subgroups for $\calP$, then the same conclusion holds.    %
\end{lem}

\begin{pf}
  Suppose otherwise, that there is a fundamental segment $\tau$ of a hyperbolic element $h$, or a path between two non-free vertices,   for which the consecutive images by $\dot f$ remain of  bounded reduced length.

Assume first that angles in the paths $\dot f^{n}(\tau)$  remain bounded, as $n$ goes to infinity.   

In this case, all paths $\dot f^{n}(\tau)$  are $\Theta$-straight   for sufficiently large $\Theta$, and bounded in length.  Therefore there are finitely possible $\dot f^{n}(\tau)$ up to the action of $G$. Let $n_2>n_1>0$ and $g\in G$ such that $\dot f^{n_2}(\tau)=g \dot f^{n_1}(\tau)$. If $\tau$ is the fundamental segment of   $\gamma$, then $g\gamma g^{-1} = \phi^{n_2-n_1}(\gamma)$, and this contradicts the fact that $\phi$ is atoroidal.
If $\tau$ is a path between two non-free vertices, it contradicts the absence of twinned subgroups.

Assume now that the angles  in the paths $\dot f^{n}(\tau)$   are unbounded.

We first treat the case of the following lemma, that we will re-use later.

\begin{lem} \label{lem;tech}  Assume that $\tau$ is a path in $\dot T$, and that the sequence of paths $\dot f^n(\tau)$ (after reduction) remains bounded in length.  Denote by $v_1(n), \dots, v_{\ell_n} (n)$ be  the consecutive non-free vertices on $\dot f^n(\tau)$, different from its end points.
   
  Assume that we may extract a subsequence   $\dot f^{n_k}(\tau)$   so that the angle at two vertices $v_{i(n_k)}(n_k), v_{j(n_k)}(n_k)$, for $i(n_k)<j(n_k)$,  is tending to infinity,    then $\phi$ has twinned subgroups.

\end{lem}
\begin{pf}
Note that $\ell_n$ is bounded by assumption. 
First   we may choose the subsequence so that $|i(n_k)-j(n_k)|$ is constant minimal.

For arbitrary $\Theta$ and $m>0$,   there is $k_0$ such that if $k>k_0$, then the angles at  $v_{i(n_k)}(n_k), v_{j(n_k)}(n_k)$    are so large that by Lemma \ref{lem-finitely-many-theta-equivalent-paths},   the map $\dot f$ applied $m'<m$ consecutive times to  the path  $\dot f^{n_k}(\tau)$ contains after reduction, the vertices $\dot f^{m'}(v_{i(n_k)}(n_k)), \dot f^{m'}(v_{j(n_k)}(n_k))$, in that order, and the angle there is still larger than $\Theta$. However, by minimality of the extracted subsequence, we may assume that the angles at the vertices between these two are bounded by $\Theta$. Thus,   the paths between $\dot f^{m'}(v_{i(n_k)}(n_k))$ and $ \dot f^{m'}(v_{j(n_k)}(n_k))$ are $\Theta$-straight, of bounded length (as subpaths of  $\dot f^n(\tau)$)  and live in a finite set. If $m$ is larger than the cardinality of this finite set, we see that we two of them are $G$-translates of one another, and this produces twinned subgroups.    
\end{pf}

Let us come back to the proof of Lemma \ref{e-g-fs-L}.  The lemma allows to exclude the case described in its assumption.
We thus now assume that there is a bound $\Theta_0$ so that all angles except perhaps one, in $\dot f^n(\tau)$, are smaller than $\Theta_0$. On the other hand, we assumed that angles were not bounded, so we may extract a subsequence and find $\dot f^{n_k}(\tau)$  so that the angle at two vertices $v_{i(n_k)}(n_k)$ goes to infinity. We may assume that $i(n_k)$ is constant, and minimal. In that way, if the initial point of $\tau$ is a non-free vertex, we obtain the same contradiction as in the previous paragraph. So we now assume that the initial point of $\tau$ is a free vertex, which means that $\tau$ is the fundamental segment of some hyperbolic element $\gamma$. Consider then $\gamma^2$. A fundamental segment  for this element consists of the concatenation of $\tau$ and $\gamma\tau$, which is a reduced concatenation, since $\tau$ is in the axis of $\gamma$. Now, images by  $\dot f$ also fail to grow, and have a pair of angles going to infinity. The previous case applies, and leads to the desired contradiction.

\end{pf}

For a path $\alpha$ in the tree $T$, $[\alpha]$ denotes its reduction in $\dot T$, a reduced path in $\dot T$ with same end points, and no collapsible subsegment.

For the following  lemma, one can also refer to \cite[Prop. 3.12]{Hor}.

\begin{lem}[Bounded cancellation lemma]

  Let $\phi\in Aut(G,\calG)$ be an automorphism of $G$, $T\in \calT_\calG$,  $f:T\rightarrow T$ be piecewise linear representing $\phi$.  Let $\calP$ be a $\phi$-invariant collection that is hyperbolically coning-off $\calG$, and $\dot f: \dot T\to \dot T$ the induced map on the cone-off. Then exist a constant $C_f$, depending only on $f$ and $\calP$, such that for any path $\rho\subset T$ obtained by concatenating two  paths $\alpha,\beta$ without cancellation, we have
$$l_{\dot T}([\dot f([\rho])])\geq l_{\dot T}([\dot f([\alpha])])+l_{\dot T}([\dot f([\beta])])-C_f$$
\end {lem}

\begin{pf}
Observe that 
 $\dot f: \dot T\to \dot T$ is a quasi-isometry (see for instance \cite[Prop. 3.2]{SisProj}). Thus there is a distance bound $D$ on the pairs of points of $\dot f( [\alpha] ) $   and $\dot f [\beta]$ that are sent at distance $\leq 2\delta$ (for $\delta$ the hyperbolicity constant of $\dot T$).  Consider the two paths $[\dot f([\alpha])]$   and $[\dot f([\beta])]$, and write $a,c$ and $c,b$ their respective end points in $\dot T$.  Consider points $a',b'$ on theses paths, at distance at most $\delta$ from the center of the tripod $a,b,c$ in $T$.    The distance $\Delta$, in $\dot T$ between $a$ and $b$ is, up to an error bounded by at most $2\delta$, the sum of the length of the subpaths from $a$ to $a'$ and from $b'$ to $b$.  As we initially noticed,   the distance in $\dot T$ from $a'$ to $c$ is at most $\lambda D$.  Thus the length of $[\dot f([\rho])])$ is at least $\Delta -2\lambda D$, hence at least  $  l_{\dot T}([\dot f([\alpha])])+l_{\dot T}([\dot f([\beta])])  -2\delta-2\lambda D$.
\end{pf}

Assume that $f: T\rightarrow T$ is a train track map representing $\phi\in Aut(G,\mathbb{G})$. Denote by $\lambda$ the growth rate of $f$. After taking a power of $\phi$, the map $f^m$ on the tree $T$ still is a train track map representing $\phi^m$, and its growth rate is $\lambda^m$. Thus, up to taking a power, we may assume that $\lambda >A$, for $A$ the constant of the transversality assumption of Theorem \ref{thm;rel-hyp-fully-irreducible_dot}. Note that if $\calP=\calG$ we do not need to take this power. From now on we assume $\lambda >A$.

Note that this means that for any legal segment $\rho$, in $T$ one has $l_{T} (f(\rho)) =  \lambda l_{ T} (\rho)$, and in $\dot T$, one has $l_{\dot T} ( [\dot f([\rho])] )  \geq \frac{\lambda}{A} l_{\dot T} (\rho) \geq  \frac{\lambda}{A} l_{\dot T}   ([\rho])$.

The bounded cancellation lemma allows to prove the following, by induction (we refer to the proof of Brinkmann, in \cite[Lem. 5.2]{PB}).

\begin{lem} (See  \cite[Lem. 5.2]{PB}) \label{Lemma-critical-constant-exponential-growth}
If $\beta$ is a legal path in $T$ with $\frac{\lambda}{A} l_{\dot T}([\beta])-2C_f>l_{\dot T}([\beta])$ (i.e. $l_{\dot T}([\beta])>\frac{2C_f}{\frac{\lambda}{A}-1}$), and if $\alpha, \gamma$ are paths such that the concatenation $\alpha-\beta-\gamma$ is locally injective, then there exists a constant $\nu>0$ (independent of $\beta$) such that the length of a maximal legal  segment of $\dot f^i([\alpha-\beta-\gamma])$ corresponding to $\beta$ is at least $\nu\left(\frac{\lambda}{A}\right)^il_{\dot T}(\beta)$ for all integer $i>0$.
\end{lem}

\begin{defn}(Critical constant of a train track map)
Let $G$ be a free product, with free factor system $\calG$, and $\calP$ hyperbolically coning off $\calG$. Let $\phi\in Aut(G,\mathbb{G})$ be irreducible, $f:T\rightarrow T$ be a train track map representing $\phi$, $\lambda$ be the growth rate of $f$, and $A$ the transversality constant of legal paths with respect to $\calP$.  The constant $\frac{2C_f}{\frac{\lambda}{A}-1}$ is called \emph{the critical constant} of $f$, where $C_f$ is the constant defined in the Bounded Cancellation Lemma.
\end{defn}

\subsubsection{Legal Control in Iteration}

\begin{lem}\label{lem;concatenate} (Analogue of  Lemma \cite[Lem. 2.9]{BFH})
   
  The concatenation in $T$ of two (pre-)Nielsen paths $\rho_1, \rho_2$ whose only common point is a free vertex is still (pre-)Nielsen.
\end{lem}

\begin{pf} Let $\rho_1, \rho_2$ be Nielsen paths with only one common point, end point of $\rho_1$, and starting point of $\rho_2$, which we   write as $v$.  
  By  assumption, there exist some $N, g_1, g_2$ such that $[f^N(\rho_1)]=g_1\rho_1, [f^N(\rho_2)]=g_2\rho_2$. It follows that $f^N(v) = g_1v = g_2v$. Since the stabilizer of $v$ is trivial, $g_1=g_2$.  This implies that $\rho_1-\rho_2$ is still a Nielsen path.  

  If the paths are only pre-Nielsen, then there are $n_i$ such that $f^{n_i} (\tau_i)$ are Nielsen paths, and all their images by $f$ too. Thus $f^{n_1\times n_2}  (\rho_1- \rho_2 )$  is a Nielsen path, and  $\rho_1- \rho_2$ there is pre-Nielsen.
\end{pf}

Following \cite{BFH}, the next lemma will be  proved by a similar idea, but an angle analysis is used to overcome the obstacle of the non-local finiteness of the involved tree.

\begin{lem} \label{3wgL}
Let $G$ be a group with  a free factor system $\calG$, of Scott complexity different from $(1,1), (0,2)$,  and $\calP$ hyperbolically  coning off $\calG$.    
 Let $\phi\in {\rm Aut}(G,\calG)$ be fully irreducible, without twinned subgroups in $\calP$, transversal to legal paths for $\phi$. Let  $f:T\rightarrow T$ be a train track map representing $\phi$. 

 Then for every $C>0$, there exist an exponent $M>0$, such that for any path $\rho$ in $T$, one of the three following holds:
\begin{itemize}
  \item (A) the length  in $\dot T$ of the longest legal segment of the reduction   $[f^M(\rho)]$  is greater than $C$;
  \item (B) the reduction in $T$ of $f^M(\rho)$ has strictly less illegal turns than $\rho$;
  \item (C) the reduced path $[\rho]$ is equal to a concatenation of $\gamma_1-\alpha_1-\dots-\alpha_s-\gamma_2$, where $\gamma_1, \gamma_2$ has at most 1 illegal turn with length at most $2C$, and that each $\alpha_i$ is a pre-Nielsen path with at most 1 illegal turn, and the end points of the $\alpha_i$ are all free vertices, except at most one.
\end{itemize}
\end{lem}

\begin{pf}
Assume (B) fails for all integer $M>0$, then no illegal turn becomes legal after iteration. In addition, let us assume that (A) fails as well. As $f$ is a train track map, none of the legal turns become illegal, the total number of illegal turns (and henceforth the number of legal segments) thus remains the same after iteration. Since $(A)$ fails,  each legal segment has a uniformly bounded length (in $\dot T$) after iteration, then there is an exponent $N$ such that, if $\pi_{\dot T} : \dot T \to G\backslash \dot T$ is the quotient map,
$$\pi_{\dot T}(\rho)=\pi_{\dot T}(\dot f^N(\rho))=...=\pi_{\dot T}(\dot f^{iN}(\rho))=...$$
for all $i\in \mathbb{Z}$.

We classify $\rho$ in the following cases.

Assume first that angles at every vertex in $[f^n(\rho)]$ remain bounded. The argument takes place in $T$ for this case. 

If statements (A) and (B) fail, the length of $\rho$ (in $T$) is bounded after iteration (by assumption of transversality of $\calP$ to legal paths, or alternatively, by bound on the angles at vertices of the cone-off by $\calP$). 
Hence, if angles are bounded too, there are only finitely many possible $G$-orbits of paths, for  the reduction in $T$ of  $f^n(\rho)$.        
Hence, there exist $N_0>0, n>0, g\in G$ such that, for reductions in $T$,  $f^n(f^{N_0}(\rho))=g f^{N_0}(\rho)$ (in other word, $\rho$ is pre-Nielsen). Denote by $v_{i,1}, v_{i,2}$ the starting and ending vertices of the maximal legal segment $\rho_i$.

Apply $ f^n$ to each maximal legal segment $\rho_i$ in $f^{N_0}(\rho)$. Since these segments are legal, the length of $f^n(\rho_i)$  grows by a factor $\lambda^n$.  These paths cancel at the possible illegal turns with other connected maximal legal segment(s), and reduces to $g \rho_i$ for the element $g$ above.   The legal segment between $gv_{i,1}$ and $f(v_{i,1})$ (if they are different) and the legal segment between $gv_{i,2}$ and $f(v_{i,2})$ (if they are different) are canceled. Hence there is a subsegment (which is legal) $\zeta_i$ of $\rho_i$ such that $f(\zeta_i)\subset g\rho_i$. For this reason, there is a vertex $v_i$ in each $\zeta_i$ (thus it is in $\rho_i$) such that $[f^n(f^{N_0}(v_i))]=g [f^{N_0}(v_i)]$.
We have that $\rho$ is a concatenation of $\gamma_1-\alpha_1-\dots-\alpha_s-\gamma_2$, where $\gamma_1, \gamma_2$ have at most 1 illegal turn with length at most $2C$, and that each $\alpha_i$ is a pre-Nielsen path with at most 1 illegal turn. In addition, as $\rho$ pre-Nielsen,  $\gamma_1,\gamma_2$ are also pre-Nielsen. Finally observe that the reduction $[\rho]$ is obtained by taking   a concatenation of reductions 
$[\gamma_1]-[\alpha_1]-\dots-[\alpha_s]-[\gamma_2]$ and reducing it further. Assume that this concatenation is not reduced (that one has to reduce it further), it means that for some indices $i$, $\alpha_i$ and $\alpha_{i+1}$ respectively end and start by the same edge adjacent to a vertex of the cone-off. One can then  replace $\alpha_i$  by the path that is equal to it except the last edge that is replaced by the first  edge of $\alpha_{i+1}$, and then replace  $\alpha_{i+1}$ by the path that is equal to it except with the first edge removed.  This still satisfy the conditions, and has less non-reduced points, thus after finitely many such changes, one get a reduced path, hence equality with $[\rho]$.

Assume now that the sequence  $[f^n(\rho)]$ have unbounded angles. We may apply Lemma \ref{lem;tech}, in order to reduce to the case that there exists $\Theta_0$ for which  at most one vertex of $[f^n(\rho)]$ has angle larger than $\Theta_0$. Let us extract a subsequence, so that there is $v_1$ in $[\rho]$  for which the image by $\dot f^{n_k}$ has angle larger than $\Theta_0$.

Denote the starting and ending vertex of $[\rho]$ by $v_a, v_b$ respectively, subdivide the path $[\rho]$ into two segments  $\rho_1= [v_a,v_1], \rho_2= [v_1,v_b]$.  Since $[f^{n_k}(\rho)]$ has a large angle at $\dot f^{n_k}(v_1)$,  the segments  $[f^{n_k}(\rho_1)]$ and  $[f^{n_k}(\rho_2)]$ make a reduced concatenation at $\dot f^{n_k}(v_1)$.

It follows that both segments remain short, and have no large angle. Thus $\rho_1$ and $\rho_2$ are pre-Nielsen.

By induction on the length,  we can further subdivide $[\rho_1]$ and $[\rho_2]$ such that
$[\rho]$ is a concatenation of $\gamma_1-\alpha_1-\dots-\alpha_s-\gamma_2$, where $\gamma_1, \gamma_2$ has at most 1 illegal turn with length at most $2C$, and that each $\alpha_i$ is a pre-Nielsen path with at most 1 illegal turn.

To see that at most one of the end points of $\alpha_i$ can be non-free, assume that two of them are, we thus have $\tau = \alpha_i-\dots-\alpha_j$, for $i\leq j$, which is a concatenation of   pre-Nielsen paths with only free vertices as intermediate subdivision points.  Their concatenation is still pre-Nielsen,  by Lemma \ref{lem;concatenate}. However the end points of $\tau$ are non-free vertices, therefore by Lemma \ref{e-g-fs-L} the iterates of $f$ on $\tau$ are eventually arbitrarily long. This contradicts that they are periodic. Thus,  at most one of the end points of $\alpha_i$ can be non-free. 

In conclusion, statement (C) of the lemma holds.
\end{pf}

\begin{lem}\label{lem;no_concatenation}
Let $G$ be a group with a free factor  system $\calG$ of Scott complexity is different from $(1,1), (0,2)$, and $\calP$ as before, coning-off $\calG$, 
and  $\phi\in {\rm Aut}(G,\calG)$ fully irreducible, preserving $\phi$.

Assume  that $\phi$ is atoroidal for $\calP$, and without twinned subgroups for $\calP$, and transversal to legal paths for $\phi$.  %

If $f:T\to T$ is a train track representative for $\phi$ on $T$. Then there exist a constant $M_0$ such that  for  any reduced (in $\dot T$)  concatenation of $M_1$  pre-Nielsen paths in $\dot T$ whose end points are free vertices,   one has $M_1<M_0$. 
\end{lem}
\begin{pf}   Recall that by Lemma \ref{lem;concatenate},  the concatenation of  such pre-Nielsen paths is still a pre-Nielsen path.

In order to prove the lemma, it suffices to prove that after concatenating sufficiently many pre-Nielsen paths in $\dot T$, in such a way that the concatenation is reduced in $\dot T$,  any resulting path contains a subconcatenation that will grow eventually after iteration of $\dot f$. This lead to a contradiction, since   this subconcatenation is still a concatenation of Nielsen paths, hence the length remains the same after iteration).   

Let us write $[\tau]= [\rho_1]-[\rho_2]-\dots-[\rho_{n_0}]$ our concatenation in $\dot T$, and $e_i$ the initial oreiented edge of  $\rho_i$.
The graph $G\backslash \dot T$ is finite, thus if $n_0$ is larger than its number of oriented edges, there is
$i<j$ such that the path $\tau'= [\rho_i]-[\rho_{i+1}]-\dots-[\rho_{j-1}]$ starts by the same free vertex and oriented edge than $\rho_j$. Thus, there is an hyperbolic element $g\in G$ that sends the initial vertex of $\tau'$ to the initial vertex of $\rho_j$, and that has $\tau'$ as fundamental segment ($\tau'-g\tau' $ is reduced). Since $\tau'$ is a Nielsen path,  it contradicts  Lemma \ref{e-g-fs-L}.

\end{pf}

As an application of Lemmas \ref{3wgL}, \ref{lem;no_concatenation}  and \ref{Lemma-critical-constant-exponential-growth} we have:

\begin{lem} (Analogue of \cite[Lem. 2.10]{BFH})\label{sbpL}

Let $\phi\in Aut(G,\calG)$ and assume that   $f:T\rightarrow T$ and $f':T'\rightarrow T'$ are train track maps representing $\phi$ and $\phi^{-1}$ respectively, that     satisfy the conclusion of Lemma \ref{3wgL}, Lemma \ref{lem;no_concatenation}, and Lemma \ref{Lemma-critical-constant-exponential-growth}, and such that $\lambda/A >1$. %
 Let $\alpha: T\rightarrow T'$ and $\alpha': T'\rightarrow T$ be $G$-equivariant, and Lipschitz.

 Then for any $C>0$ there exist exponent $N>0$ and $L_0>0$, such that if $\rho$ is a path whose reduction in $\dot T$ is of length $\geq L_0$, and $\rho'$ is a reduction  in $\dot T'$ of
 $\alpha(\rho)$, then either $[f^N(\rho)]$ or $[f'^N(\rho')]$ contains a legal segment of length greater than $C$.
\end{lem}

\begin{pf}
 We first gather a few constants before starting the proof.   By the above lemma  \ref{lem;no_concatenation}, there exist a constant $M_0$ such that it is impossible to concatenate more than $M_0$ pre-Nielsen paths with end points being free vertices, and the concatenation being reduced in
 $\dot T$. 

 Fix $C>0$ such that it is larger than the critical constant for both $\dot f$ and $\dot f'$ (see Lemma \ref{Lemma-critical-constant-exponential-growth}).

 Let $M$ be the greater one of the integers according to Lemma \ref {3wgL} when we apply it to $f, C$ and $f', C$.

 We may assume, after global rescaling, that the length of the shortest edge in $T$ and in $T'$ is larger than $1$.

 Fix $P > (2M_0+2)\times C + 2\lambda (C_f +C_{f'})$, and  $Q$ such that $\frac{P-1}{P}<Q<1$.


We argue by contradiction, assuming that for all $N$ and all $L_0$ there is a path $\rho$ that is a counterexample to the statement,
a path   whose reduction in $\dot T$ is longer than $L_0$, and such that for all $N\leq N_0$, neither the reduction  $[f^N(\rho)]$, nor $[f'^N(\rho')]$   contain a legal segment of length greater than $C$. We want to show that, for sufficiently large $L_0$ and $N_0$, this leads to a contradiction. We thus fix $L_0 >10\times C\times (2M_0+2)$  and $N_0>0$, and extra conditions on their size (to reach a contradiction) will appear at the end.


We may represent $[\rho]$   %
as a concatenation of segments $\beta\subset [\rho]$ such that $\dot f^M(\partial \beta)\subset [f^M(\rho)]$, each has at least $2M_0+2$ illegal turns, and  whose length is at most $P$.  

Then (A) of Lemma \ref {3wgL} fails for such a subpath $[\beta]$. In addition, the case (C) of Lemma \ref {3wgL} also fails: there are at most $2M_0$ subpath $\alpha_i$ in the notations of $C$, and one of them has to be longer than $2C$, implying that there is a legal subsegment of length $C$, which we excluded. 

Then (B) of Lemma \ref {3wgL} must hold for $[\beta]$.

The constant $P$ is also an upper bound to the number of illegal turns in each of these segments. 
Recall that  $\frac{P-1}{P}<Q<1$.
For a path $\tau$ we denote by $NIT(\tau)$ the number of illegal turns in $\tau$.

Then for these segments 
 $\beta$   we have
$$\frac{NIT([\dot f^M(\beta)])}{NIT(\beta)}\leq Q$$

This inequality is thus also true for any concatenation of such segments in $\rho$, hence for any sufficiently long subsegment of $\rho$.

We do the same construction to $f^M(\rho), f^{2M}(\rho),\dots$, as long as $uM \leq N_0$.   Then for  all such $u$, and all  segments $\beta$ in $[\rho]$ whose lenght is sufficiently large (say, for notation,  larger than $L(u)$), 
$$\frac{NIT([\dot f^{uM}(\beta)])}{NIT(\beta)}\leq Q^u$$

Of course for such a $\beta$ to exist, one must have $L_0\geq L(u)$. 

Since we require that, for all $N\leq N_0$,  any legal segment in each $[f^N(\rho)]$ is bounded by $C$, and it is obviously not less than the length of shortest edge (which has length larger than $1/2$),  one has, for $uM<N_0$, and for all such $\beta$ larger than $L(u)$,    
$$\frac{l_{\dot T}([\dot f^{uM}(\beta)])}{l_{\dot T}(\beta)}\leq 2(C+1) Q^u \eqno{(1)}$$

Apply the same discussion to $[\alpha f^{uM}(\rho)]$ as we did to $\rho$, and consider $f'$ instead of $f$. We then  have, for $\beta$ larger than some constant $L'(u)$, 
$$\frac{l_{\dot T'}([\dot f'^{uM}\alpha \dot f^{uM}(\beta)])}{l_{\dot T'}([\alpha \dot f^{uM}(\beta)])}\leq 2(C+1)Q^s \eqno{(2)}$$

For sake of notations, we will redefine $L(u)$ to be at least $L'(u)$.

Notice that $\dot f'^{uM}\alpha \dot f^{uM}$ is conjugate to $\alpha$, hence there is some constant $\mu>1$ such that for long segments,
$$\frac{1}{\mu}\leq \frac{l_{\dot T'}([\dot f'^{uM}\alpha \dot f^{uM}(\beta)])}{l_{\dot T'}([\alpha (\beta)])}\leq \mu \eqno{(3)}$$

Multiply ($1$), ($2$) and the inverse of ($3$) we have
$$\frac{l_{\dot T}([\dot f^{uM}(\beta)])}{l_{\dot T'}([\alpha \dot \dot f^{uM}(\beta)])} \frac{l_{\dot T'}([\alpha (\beta)])}{l_{\dot T}(\beta)}\leq \mu 4(C+1)^2Q^{2u}$$

Notice that $\frac{l_{\dot T}([\dot f^{uM}(\beta)])}{l_{\dot T'}([\alpha \dot f^{uM}(\beta)])} \frac{l_{\dot T'}([\alpha (\beta)])}{l_{\dot T}(\beta)} \geq\frac{1}{\mu}Lip(\alpha)Lip(\alpha')$,
where $Lip(\alpha),Lip(\alpha')$ is the Lipschitz constant of $\alpha, \alpha'$ respectively, we have

$$\frac{1}{\mu}Lip(\alpha)Lip(\alpha')\leq \frac{l_{\dot T}([\dot f^{uM}(\beta)])}{l_{\dot T'}([\alpha \dot f^{uM}(\beta)])} \frac{l_{\dot T'}([\alpha (\beta)])}{l_{\dot T}(\beta)}\leq \mu\times  4(C+1)^2Q^{2u}$$

Therefore, setting $K=  \frac{Lip(\alpha)Lip(\alpha')}{\mu^2\times 4(C+1)^2}$, the integer $u$ is bounded by the inequality $K \leq Q^{2u} $ and let $u_m$ the maximal value satisfying this. 

We now take the length $L_0$ to be larger than $L(u_m+1)$, and $N_0 >(u_m+1) M$, in order to have a path $\beta$ in $[\rho]$ of size $L(u_m+1)$ and  thus obtain $K \leq Q^{2(u_m+1)} $ contradicting the maximality of $u_m$.

\end{pf}

\begin{defn}($C$-legality of a path)
Given a $T\in \calT_\calG$ and a constant $C$, for any reduced path $\rho\subset T$, the $C$-\emph{legality} of $\rho$ is the ratio of the sum of lengths (in $\dot T$) of legal segments in $\rho$ longer than $C$ over the total length (in $\dot T$) of $\rho$, denoted by $LEG_{\dot T}(\rho)$.
\end{defn}

\begin{lem} (See \cite[Lem. 5.6]{BFH})\label {lrl}
  Let $\phi\in Aut(G,\calG)$ be fully irreducible for $\calG$ and atoroidal for $\calP$, without twinned subgroups in $\calP$, and with Scott complexity different from $(1,1), (0,2)$.  Let $f:T\rightarrow T$, $f':T'\rightarrow T'$ be train track maps representing $\phi$ and $\phi^{-1}$ respectively. And let $\alpha: \dot T\rightarrow \dot T'$ and $\alpha': \dot T'\rightarrow \dot T$ be Lipschitz map corresponding to difference of markings. Assume that $C$ is the larger one of the critical constant of $\dot f$ and $\dot f'$. 

  There is $\epsilon>0$ and an integer $N_1>0$ such that for every nontrivial $g\in G$ hyperbolic on $\dot T$, if $\sigma$ is a fundamental segment of $g$ in $\dot T$, or a path between two non free vertices, and if   $\sigma'$ is $[\alpha(\rho)]$ in $\dot T'$,  
  then for every $N>N_1$, either $LEG_{\dot T}(\dot f^N(\sigma))\geq \epsilon$ or $LEG_{\dot T'}(\dot f'^N(\sigma'))\geq \epsilon$.
\end{lem}

\begin{pf}
By Lemma \ref {e-g-fs-L}, there is an integer $N'$, such that $l_{\dot T}(\dot f^{N'}(\sigma))>L_0$ and  $l_{\dot T'}(\dot f'^{N'}(\sigma'))>L_0$, where $L_0$ is defined according to Lemma \ref {sbpL}. And by Lemma \ref {sbpL}, there is $N_1$ such that either $[\dot f^{N_1}(\sigma)]$ or $[\dot f'^{N_1}(\sigma')]$ contains a legal segment of length greater than $C$.

Suppose the result does not hold, then there is a sequence $\{g_i\}$ with $\{\sigma_i\}$ and $\{\sigma'_i\}$ in $\dot T$ and $\dot T'$ respectively such that the legality of $LEG_{\dot T}(\dot f^{N_1}(\sigma_i))$ and of  $LEG_{\dot T'}(\dot f'^{N_1}(\sigma_i'))$   converges to 0.

Then there exists arbitrarily long segments in $\dot f^{N_1}(\sigma_i)$ and in $\dot f'^{N_1}(\sigma_i')$ (as $i$ varies) that do not contain a legal segment of length $\geq C$. Thus contradicts the Lemma \ref {sbpL}.
\end{pf}

\subsection{Relative hyperbolicity in the fully irreducible case} \label{sec;proof_fully_irred}

\begin{pro}\label{prop;proto_rel_hyp_fully_irred_trees}
  Let $G$ be a finitely generated group with a free factor system $\calG$ of Scott complexity  $(k,p)$  different from $(1,1), (0,2)$.

  Consider  $\phi\in Aut(G,\calG)$  fully irreducible. Consider $\calP$ hyperbolically coning off $\calG$, preserved by $\phi$, transverse to the legal paths of $\phi$, and such that  $\phi$ is atoroidal for $\calP$, with no twinned subgroups for $\calP$.

 Let  $f:T\to T$  be a train track map realizing   $\phi$, and $f': T'\to T'$  a train track map realising $\phi^{-1}$. Let $\dot f: \dot T\to \dot T$ and $\dot f' : \dot T' \to \dot T'$ the maps on the cone-off for $\calP$. 

 Assume that there is $A>1$ such that if $\rho$ is a legal path in $T$, the length $l_{\dot T}([\rho])$ of its reduction in $\dot T$ is of length at least $ l_{\dot T}([\rho])/A$, and similarily for $f'$.
 
 Then, for some $m\geq 1$,  the pair   $(\dot f^m, (\dot f')^m)$ is a hyperbolic pair, in the sense of Definition \ref{def;hyp_pair}.

\end{pro}

\begin{pf}
By Lemma \ref{lem;stretching_factor_Scott} the growth rate of $f$ and $f'$ are greater than some $\lambda>1$. Up to taking some power, we may assume that $\lambda >A$. 
  
   Both trees being in $\calT_\calG$, coned by the same family, there exists $\alpha: \dot T\rightarrow \dot T'$ and $\alpha': \dot T'\rightarrow \dot T$ two $G$-equivariant Lipschitz maps. For any hyperbolic element $x\in G$, let $\sigma$ be a fundamental segment of $x$ in $\dot T$,  or a path between non free vertices, let   $\sigma'$ be $[\alpha(\rho)]$, the reduction in $\dot T'$ of $\alpha(\rho)$.

Choose $C$ be the larger one of the critical constant of $\dot f$ and $\dot f'$. By Lemma \ref {lrl}, there is $\epsilon>0$ and integer $N_1>0$ (regardless of the choice of $x$ and $\sigma$) such that for every $N>N_1$, either $LEG_{\dot T,C}(\dot f^N(\sigma))\geq \epsilon$ or $LEG_{\dot T',C}(\dot f'^N(\sigma'))\geq \epsilon$.

Denote here by $S_C(\sigma)$ the set of all maximal legal segments in $\sigma$ with length longer than $C$.

We assume that $LEG_{\dot T}(\dot f^N(\sigma))\geq \epsilon$ (the other case is similar). By definition of critical length and  Lemma \ref{Lemma-critical-constant-exponential-growth}, there is $\nu>0$ such that for all $i\geq N_1$, 
$l_{\dot T}([\dot f^i(\sigma)])\geq \nu(\frac{\lambda}{A})^i l_{\dot T}(S_C(\sigma))$.
   
Moreover, we have that  $  l_{\dot T}(S_C(\sigma)) \geq \epsilon\times  l_{\dot T}(\sigma)$. It follows that   $l_{\dot T}([\dot f^i(\sigma)])\geq    (\nu\epsilon)(\frac{\lambda}{A})^i  l_{\dot T}(\sigma)    $ for all large $i$.

For sufficiently large $i$,  this establishes that $(\dot f, \dot f')$ is a hyperbolic pair. 

\end{pf}

By Remark \ref{rem;pass} and the previous proposition,  we thus have the following.
\begin{pro}\label{prop;proto_rel_hyp_fully_irred}
  Let $G$ be a finitely generated group with a free factor system $\calG$ of Scott complexity $(k,p)$   different from $(1,1), (0,2)$.

 Consider  $\phi\in Aut(G,\calG)$  fully irreducible. Consider $\calP$ hyperbolically coning off $\calG$, preserved by $\phi$, transverse to the legal paths of $\phi$, and such that  $\phi$ is atoroidal for $\calP$, with no twinned subgroups for $\calP$.

 Then $\phi$ is hyperbolic relative to $\calP$, in the sense of Definition \ref{def-rel-hyp-aut}.
\end{pro}

We may finally prove Theorem \ref{thm;rel-hyp-fully-irreducible_dot}:  it is a direct consequence of  Proposition  \ref{prop;proto_rel_hyp_fully_irred} and Gautero and Weidmann's Theorem  \ref{GW-comination-theorem-full}.

Alternatively, one may use Mj and Reeves's Theorem \ref{thm;MjR} (see also \cite[Thm. 2.20]{Gau}) instead of  \ref{GW-comination-theorem-full}. We briefly indicate how, without entering in the detail of the definitions from \cite{MjR}.

In our case, the finite graph of groups is a loop of groups, whose vertex and edge groups are $G$, and attaching maps are given by the identity and $\phi$. All properties from $1$ to $3$ are obviously satisfied. The induced tree of spaces, is just a bi-infinite line of spaces, and at each integer, the space is the cone-off tree $\dot T$ associated to the train track, with attachments with next cone-off tree, and the previous  being given, respectively,  by $\dot f$ and $\dot f_{-1}$, a continuous map from $\dot T$ to itself realizing $\phi^{-1}$. The hallways flare condition \cite[Def. 3.3 -- 3.5]{MjR} is  the expansion (under the power of $\dot f$, or of $\dot f_{-1}$ )  of sufficiently long paths in $\dot T$, which is ensured by the expansion of the paths corresponding to fundamental segments of hyperbolic elements, from Proposition \ref{prop;proto_rel_hyp_fully_irred_trees}, and the cone-bounded hallways strictly flare condition  is the expansion (under the power of $\dot f$, or of $\dot f_{-1}$ ) of the paths between two non-free vertices, from     Proposition \ref{prop;proto_rel_hyp_fully_irred_trees}. 
The collection $\calC$ corresponds, in \cite{MjR},  to our mapping torus of the free factor system.

Let us finally make the comment that, using the theorem of Gautero-Weidmann does not require the control of the expansion of the paths between non-free vertices, and therefore can be done without the assumption of absence of twinned subgroups. We do not know whether this makes a difference.

\subsection{Relative hyperbolicity in the reducible case}

Let $G$ be a group with a free factor system $\calG$, and consider an automorphism $\phi\in {\rm Aut} (G, \calG)$,  that is atoroidal, but possibly  reducible with respect to  $\calG$.

\begin{thm}\label{thm;no_twin} Let $G$ be a finitely generated group, and $\calG$ be a free factor system. Let $\phi \in {\rm Aut} (G,\calG)$ be atoroidal for $\calG$.  Assume that there is no pair of twinned subgroups in $G$ for $\calG$ and $\phi$.

Then $G\rtimes_\phi \mathbb{Z}$ is relatively hyperbolic with respect to the mapping torus of $\calG$. 
\end{thm}

\begin{pf} Let  $\calG_{{\rm fi}}$ be  the free factor system  of $G$ provided by Lemma \ref{lem-existence-fully-irre-ffs}: it is preserved by some power $\phi^{m}$ of $\phi$, and $\phi^{m}$ is fully irreducible with respect to $\calG_{{\rm fi}}$.

We claim that  for some $m'$   $G\rtimes_{\phi^{m'}} \mathbb{Z}$ is relatively hyperbolic with respect to the mapping torus of $\calG$. We prove this     by discussing according to the Scott complexity of  $G$ for  $\calG_{{\rm fi}}$.

If the Scott complexity of $G$ for $\calG_{{\rm fi}}$ is $(0,2)$, then  $\calG_{{\rm fi}}= \{[A],[B]\}$ and $G=A*B$.   After composing by a conjugation, we may assume that  automorphism $\phi^m$  preserves $A$, and then it sends $B$ on a conjugate of $B$ by an element of $A$. By composing by another conjugation we may assume that it preserves both $A$ and $B$. Thus, $G\rtimes_{\phi^m} \mathbb{Z}$ is isomorphic to $(A\rtimes_{\phi^m|_A}  C )*_{C} (B\rtimes_{\phi^m|_B} C)$, where $C$ is the infinite cyclic group generated by the element associated to $\phi^m$.   As we have proved,   both factors are relatively hyperbolic, with respect to  mapping tori of the free factor systems induced by $\calG$. Since $\phi$ (hence $\phi^m$) has no twinned subgroups in $\calG$, the group $C$ is non-parabolic in at least one of the factors. The combination theorem  \cite[Thm. 0.1, case 3]{Dah03}  can thus be applied, and  the group $G$ is hyperbolic relative to the mapping torus of the union of the free factor systems of $A$ and $B$ induced by $G$, hence to the mapping torus of $\calG$ for $\phi^m$. We have the claim for this case.

If  the Scott complexity of $G$ for $\calG_{{\rm fi}}$ is $(1,1)$, then  $\calG_{{\rm fi}}= \{[A]\}$ and $G= A* Z$  for  a subgroup $Z$, infinite cyclic, generated by $z$.  
 We may assume, after composing by a conjugation, that $\phi$ preserves $A$.  Also, since $\phi$ is an automorphism,  there exists $g\in A$ such that $\phi^2(z)$ is $z g$ (the square ensures that the exponent of $z$ is $+1$ instead of $-1$). 
One can thus express 
$G\rtimes_{\phi^2} \langle t \rangle  $ as isomorphic to $  \left( A\rtimes_{\phi^2|_A}  \langle t \rangle \right)  *_{\langle t \rangle,   \langle tg^{-1} \rangle  }$.
 As we have proved, $A\rtimes_{\phi|_A}  \langle t \rangle $ is  relatively hyperbolic with respect to  mapping torus of the free factor system induced by $\calG$. 

We claim that either $t$ or $tg^{-1}$ is not parabolic in $A\rtimes_{\phi|_A}  \langle t \rangle$. If both are parabolic, there are two free factors  of $A$ in $\calG$, say $H, K$, such that $H$ is normalized by $t$ and $K$ by $tg^{-1}$.  Since $z^{-1} t z = t g^{-1}$, we have that $tg^{-1}$ normalises both $K$ and $z^{-1} H z$. By the absence of twinned groups for $\phi^2$, it follows that $K= z^{-1} H z$. But considering normal forms for the free product $G= A*Z$, we have that  $z^{-1} H z$ is not a subgroup of $A$.  This proves that either $t$ or $tg^{-1}$ is not parabolic in $A\rtimes_{\phi|_A}  \langle t \rangle$.  We may therefore apply the combination theorem  \cite[Thm. 0.1, case 4]{Dah03} to obtain that the HNN extension  $  \left( A\rtimes_{\phi^2|_A}  \langle t \rangle \right)  *_{\langle t \rangle,   \langle tg^{-1} \rangle  }$ is relatively hyperbolic with respect to the conjugates of the parabolic subgroups of  $\left( A\rtimes_{\phi^2|_A}  \langle t \rangle \right) $. In other words $G\rtimes_{\phi^2} \mathbb{Z}$ is hyperbolic relative to the mapping torus of $\calG$ for $\phi^2$.

If the Scott complexity of $G$ for  $\calG_{{\rm fi}}$ is different from $(1,1)$ and $(0,2)$, then $G\rtimes_{\phi^{m}} \mathbb{Z}$ is relatively hyperbolic with respect to the mapping torus of $\calG_{{\rm fi}}$       by  Theorem \ref{thm;rel-hyp-fully-irreducible}.  

We argue then by induction on the Scott complexity of $(G,\calG)$. The lowest complexities $(0,2)$ and $(1,1)$ have just  been treated.    Consider $H$ such that $[H] \in \calG_{{\rm fi}}$, and let $g_H$ such that ${\rm ad}_{g_H} \circ \phi^{m}$ preserves $H$.  Let $\calH$ be the free factor system of $H$ induced by $\calG$.    
One can easily check that the automorphism ${\rm ad}_{g_H} \circ \phi^{m}$ of $H$ is atoroidal  and has no twinned subgroups  for $\calH$. 
Also, the Scott complexity of $H$ for the free factor system $\calH$ is strictly lower than that of $\calG$ for $G$, by Lemma \ref{kurosh-scott}. By induction hypothesis, for each element $H$ of $\calG_{{\rm fi}}$, its mapping torus by $({\rm ad}_{g_H} \circ \phi^{m_H})|_H$ is hyperbolic relative to the mapping torus of the free factor system induced by $\calG$.

We can conclude about the claim for $G$,    by the telescoping argument of Osin \cite[Thm 2.40]{Osi}:  
  the group     $G\rtimes_{\phi^{m'}} \mathbb{Z}$ is relatively hyperbolic with respect to the mapping torus of $\calG$, for $m'$ the product of $m$ and of the exponents $m_H$ associated to the groups $H$ in $\calG_{{\rm fi}}$.

We thus have proved the claim for all cases.

Finally, by   Proposition \ref{prop;phimtophi}  we get that the  group $G\rtimes_{\phi} \mathbb{Z}$ is relatively hyperbolic with respect to the mapping torus of $\calG$. %

\end{pf}

\section{Applications}
\subsection{The central element condition}\label{sec;31}

\begin{pro} \label{prop;no_twin} Assume that $\phi\in {\rm Aut}(G,\calG)$ satisfies the \emph{central element condition}:
for all $[H]\in \calG$, there exists $g\in G$ such that ${\rm ad}_g\circ \phi|_H$ is an automorphism of $H$, there is a non-trivial central element in $H\rtimes_{{\rm ad}_g\circ \phi|_H} \mathbb{Z}$.

Then, if $\phi$ is atoroidal, then it has no twinned subgroups for $\calG$.   
\end{pro}

\begin{pf}  Assume it has a pair of twinned subgroups. There exists $g\in G$, and $m\geq 1$ such that $A,B$, both fixed by ${\rm ad}_g\circ \phi^m$. By the central element condition, there exists $a, b$ respectively  in $A, B$,   non-trivial, so that $ {\rm ad}_g\circ \phi^m(a) = a$ and ${\rm ad}_g\circ   \phi^m(b) = b$. The product $ab$ is in $Hyp(\calG)$ but is fixed by ${\rm ad}_g\circ \phi^m$. This contradicts atoroidality of $\phi$. 
\end{pf}

In the next statement, we say that an automorphism $\phi\in {\rm Aut}(G,\calG)$ is \emph{toral}, if  for each $H$ such that $[H]\in \calG$, there exists $g\in G$ such that ${\rm ad}_g\circ \phi|_H$ is the identity on $H$.
The three corollaries are consequences of Proposition \ref{prop;no_twin} and Theorem \ref{thm;no_twin}.

\begin{cor}\label{coro;toral}
Assume that $G$ is finitely generated, and that $\calG$ consists of torsion free abelian groups. 
If $\phi$ is atoroidal, and toral, then the group $G\rtimes_\phi \mathbb{Z}$ is toral relatively hyperbolic.  
\end{cor}

\begin{cor}
Assume that $G$ is torsion free, and that $\calG$ consists of nilpotent groups and that $\phi\in {\rm Aut}(G,\calG)$  is such that for each $[H]\in \calG$, there exists $g\in G$ such that ${\rm ad}_g\circ \phi|_H$ is the identity on $H$. Then if $\phi$ is atoroidal, the group $G\rtimes_\phi \mathbb{Z}$ is  relatively hyperbolic with nilpotent parabolic subgroups.  
\end{cor}

\begin{cor}
Assume that $G$ is torsion free, and that $\calG$ consists of abelian groups and that $\phi\in {\rm Aut}(G,\calG)$  is such that for each $[H]\in \calG$, there exists $g\in G$ such that ${\rm ad}_g\circ \phi|_H$ is unipotent on $H$ (seen as $\mathbb{Z}$-module). Then if $\phi$ is atoroidal, the group $G\rtimes_\phi \mathbb{Z}$ is  relatively hyperbolic with nilpotent parabolic subgroups.  
\end{cor}

\subsection{The case of free groups, and a theorem of Gautero-Lustig}

%

\begin{thm}(Gautero-Lustig, and Ghosh)\label{thm;GLG} \cite{GL} \cite{Gho}

  If $\phi$ is an automorphism of a free group $F$ with at least one exponentially growing element, the semidirect product $F\rtimes_\phi \mathbb{Z}$ is relatively hyperbolic with respect to the mapping torus of the collection of maximal polynomially growing subgroups.  
\end{thm}  

\begin{pf}
We now assume that $\phi$ is an automorphism of a finitely generated free group $G$. We argue by induction on the rank of $G$. 

 Let $\calG$ be a maximal invariant free factor system provided by Lemma \ref{lem-existence-fully-irre-ffs}, and passing to a power we may assume that $\phi^m\in {\rm Aut} (G,\calG)$, and thus is fully irreducible.  

 By Corollary \ref{coro;23},  there is $m_1$ such that $G\rtimes_{\phi^{m_1}} \mathbb{Z}$ is hyperbolic relatively to the  mapping tori of the maximal polynomially growing subgroups of $\phi$ with respect to $T\in \calT_\calG$. Note that in the case of Scott complexity $(1,1)$ and $(0,2)$, this holds trivially because the maximal polynomially growing subgroups of $\phi$ with respect to $T\in \calT_\calG$ is $G$ itself.

Consider then $P$ a maximal polynomially growing subgroup in $G$ for $T\in \calT_\calG$.  
By Proposition \ref{prop;Levitt_four_points} %
its rank is strictly less than the rank of $G$, or perhaps, it is $G$ itself. 
Consider $T_P$ its minimal subtree in $T$, it provides a decomposition of finite graph of groups, with trivial edge groups, and with vertex groups in the collection $\calG$. The Scott complexity of this decomposition is necessarily $(0,2)$ or $(1,1)$, since otherwise, some element in $P$ would be growing exponentially fast by iterations of $\phi^{m_1}$.  Thus it is simply either a free product of two groups $G_a, G_b$  in $\calG$, or a group $G_a$  in $\calG$ free product with $\mathbb{Z}$, and this last free factor is growing polynomially for $\phi$ (for the word metric of conjugacy classes). 

In the first case, the suspension of $P$ by  $\phi^{2m_1}$  (possibly precomposed by somme inner automorphism) is isomorphic to the amalgam over some maximal infinite cyclic subgroup of  the suspensions of $G_a$ and of $G_b$.  By induction hypothesis, both suspensions of $G_a$ and $G_b$ are relatively hyperbolic with respect to suspension of polynomially growing subgroups. Writing $t$ as the  generator of the $\mathbb{Z}$ factor, over which the amalgamation is performed, if $t$ is non-parabolic in at least one of these two suspensions, the combination theorem  \cite[Thm 0.1]{Dah03}  ensures that the suspension of $P$ is relatively hyperbolic with respect to the conjugates of the parabolic subgroups in both $G_a\rtimes \mathbb{Z}$ and $G_b\rtimes \mathbb{Z}$, as wanted. If  it is parabolic in both, then $G_a$ and $G_b$ possess  polynomially growing subgroups $P_a$ and $P_b$ normalized by $t$. The group $\langle P_a, P_b\rangle$ is then polynomially growing in $P$ for the automorphism, and the combination theorem guaranties that $P\rtimes \mathbb{Z}$, which is isomorphic to 
\[ \Big(  \left( G_a\rtimes \mathbb{Z} \right) *_{P_a\rtimes \mathbb{Z} } \left(  \langle P_a, P_b\rangle\rtimes \mathbb{Z}   \right)  \Big) *_{ \langle P_a, P_b\rangle\rtimes \mathbb{Z}   } \Big( \left( G_b\rtimes \mathbb{Z} \right) *_{P_a\rtimes \mathbb{Z} } \left(  \langle P_a, P_b\rangle\rtimes \mathbb{Z}   \right)\Big) \] is indeed relatively hyperbolic with respect to suspensions of polynomially growing subgroups.

In the second case,  the suspension of $P$ by  $\phi^{m_1}$ is isomorphic to the HNN of the suspension of $G_a$ over two maximal infinite cyclic subgroups, with the stable letter being an element of $P$, growing polynomially for $\phi$ (for the word metric of conjugacy classes).

 For each group $H$ whose conjugacy class is in $\calG$, there is $h\in G$ such that ${\rm ad}_h \circ \phi^{m_H}$ preserves $H$. 
By induction hypothesis, since $H$ has smaller rank than $G$, the group   $H \rtimes_{{ \rm ad}_h \circ \phi^{m_H}} \mathbb{Z}$ is relatively hyperbolic with respect to its polynomially growing subgroups for ${ \rm ad}_h \circ \phi^{m_h}$ (for its Cayley tree).

If the edge group is not parabolic in at least one of its embedding in its adjacent vertex group, the  combination theorem \cite[Thm 0.1]{Dah03} directly ensures that, for $m_2=2m_1m_{G_a}m_{G_b}$,  the suspension of $P$  by $\phi^{m_2}$  is relatively hyperbolic with respect to the suspensions of the maximal polynomially growing subgroups in $P$.   

If it is parabolic in both embeddings, then first perform the amalgam, over the edge group,  of $G_a$ with the maximal parabolic group containing the edge group of $G_b$ (or with an abstract copy of the  maximal parabolic group containing the second image of the edge group of $G_a$) thus enlarging the parabolic subgroup of $G_a$, but keeping the relatively hyperbolic structure. Then, perform the amalgam with $G_b$ over a maximal parabolic group in $G_b$, or the HNN over the maximal parabolic group in $G_a$, using again the combination theorem \cite[Thm 0.1]{Dah03} to obtain that 
 the suspension of $P$  by $\phi^{m_2}$  is relatively hyperbolic with respect to the suspensions of the maximal polynomially growing subgroups in $P$.

Therefore, in all cases,  the suspension of $P$  by $\phi^{m_2}$  is relatively hyperbolic with respect to the suspensions of the maximal polynomially growing subgroups in $P$. 

By the telescopic argument  of Osin \cite[Thm 2.40]{Osi}, we conclude that for some $m'$,  $G\rtimes_{\phi^{m'}} \mathbb{Z}$ is relatively hyperbolic with respect to suspensions of maximal polynomially growing subgroups.

By the reduction of Proposition \ref{prop;phimtophi} 
we obtain the result. 

\end{pf}

\subsection{Application to the conjugacy problem}

In this section, we consider the conjugacy problem for certain outer automorphisms of a group $G$  that is a  free products of non-cyclic  abelian groups.  We present an algorithm that will decide whether two given toral atoroidal automorphisms of free products of abelian groups are conjugate in $Out(G)$. Hence we will prove the following.

\begin{thm} \label{thm;CP_toral} Let $G$ be a finitely generated free product of non-cyclic free abelian groups,  $G=A_1*\dots*A_p$. Denote by $\calA$ the free factor system $\{[A_i], i=1,\dots p\}$
 
  There is an algorithm that,   given $\phi_1, \phi_2$, two automorphisms of $(G, \calA)$ that are atoroidal, and toral, determines whether they are conjugate in $Out(G,\calA)$. 
\end{thm}  

Note that we assume that the Scott complexity of $(G, \calA)$  in the statement is $(0,p)$.     %

\subsubsection{Reduction to an orbit problem for the group ${\rm Out}(  G\rtimes_{\phi_2} \langle t_2\rangle )$}

We will write $\overline{w}$ for the image of an object $w \in G\rtimes \mathbb{Z}$ in the  abelianisation of  $G\rtimes \mathbb{Z}$.   We consider $G_1, G_2$, two copies of $G$, (canonically isomorphic to $G$).  

\begin{lem} \label{orbit} $\phi_1, \phi_2 \in {\rm Aut}(G)$ are conjugated in $Out(G)$ if and only if the following holds:
  \begin{itemize}
     \item there exists an isomorphism $\alpha: G_1\rtimes_{\phi_1} \langle t_1\rangle \to G_2\rtimes_{\phi_2} \langle t_2\rangle$
     \item and there exists an automorphism $\psi$ of $G\rtimes_{\phi_2} \langle t_2\rangle$ such that, in the abelianisation,  if    $\bar G_2$, $\bar t_2$, $ \overline{\psi(\alpha (t_1) )} $ and $\overline{\psi(\alpha (G_1) )}$
       are the images of $G_2$, $t_2$,  $\psi(\alpha (t_1) )$  and    $\psi(\alpha (G_1))$ in the abelianisation of    $G_2\rtimes_{\phi_2} \langle t_2\rangle$, then
       \begin{equation} \label{inclusion}
        \bar G_2 \subset   (\overline{\psi(\alpha (G_1) )} )  \quad \hbox{ and }  \quad  (\overline{\psi(\alpha (t_1) )} ) \in   \bar t_2  \bar G_2.
         \end{equation}
  \end{itemize}
  
\end{lem}   

\begin{pf}
By \cite[Lem. 3.1]{DTAMS},  $\phi_1, \phi_2 \in {\rm Aut}(G)$ are conjugated in $Out(G)$ if and only if there is an isomorphism $\iota$  between $G_1\rtimes_{\phi_1} \langle t_1\rangle$ and  $G_2\rtimes_{\phi_2} \langle t_2\rangle$   such that at the level of the abelianisations,  $\overline{\iota ( G_1)}$ contains (hence is equal to)  $\bar G_2$, and  $\overline{\iota( t_1)}$ is  in $\bar t_2  \bar G_2$.    The Lemma follows.
\end{pf}

\begin{lem} The first point of the characterisation in the previous lemma is decidable: there is an algorithm that given $G, \calG, \phi_1, \phi_2$ indicates whether there exists such an $\alpha$. 
  \end{lem}
\begin{pf}
  If  $\calG$ is a free factor system of $G$ consisting of torsion free abelian groups, and  $\phi_1, \phi_2 \in {\rm Aut}(G,\calG)$  are atoroidal for $\calG$ and toral in the sense of  Corollary \ref{coro;toral}, then $G_1\rtimes_{\phi_1} \langle t_1\rangle$ and  $G_2\rtimes_{\phi_2} \langle t_2\rangle$ are relatively hyperbolic with torsion free abelian parabolic subgroups, in other words they are toral relatively hyperbolic. By \cite[Thm. D]{DGr} we have an algorithm determining whether they are abstractly   isomorphic, and
    if they are, we may find an isomorphism $\alpha :   G_1\rtimes_{\phi_1} \langle t_1\rangle \to  G_2\rtimes_{\phi_2} \langle t_2\rangle$.  
\end{pf}

The problem is now, assuming that such an isomorphism $\alpha$ is given,  to decide the second point of the characterisation, which we call our orbit problem in  ${\rm Aut}(  G_2\rtimes_{\phi_2} \langle t_2\rangle )$.

Observe that the properties involved in the orbit problem only depend on the class of $\psi$ in  
 ${\rm Out}(  G_2\rtimes_{\phi_2} \langle t_2\rangle )$.

We will need this decomposition lemma.   %

\begin{lem}
 If $G= A_1*\dots* A_p $ where $A_i$ are abelian, and if $\phi_2$ is an automorphism of $(G,\calG)$, the image $\bar G$ of $G$ in the abelianisation of    $G\rtimes_{\phi_2} \langle t_2\rangle$ splits as a direct product $B_1\times B_2\times \dots \times B_p$, where each $B_i$ is the image of $A_i$.
\end{lem}
 
\begin{pf}
The abelianisation of $G$  is  $G_{ab} = B_1\times \dots\times  B_p$, where $B_i$ is canonically isomorphic to $A_i$.   Let $\bar \phi_2$ the induced automorphism of $G_{ab}$.  Since $\phi_2$ preserves the conjugacy class of $A_i$,   the automorphism $\bar \phi_2$ preserves each $B_i$ in $G_{ab}$.  Write $C_i$ for the abelianisation of $B_i\rtimes_{\bar \phi_2 |_{B_i}}  \mathbb{Z}$, and $c_i\in C_i$ the element associated to the factor $\mathbb{Z}$. 
The abelianisation of  $G\rtimes_{\phi_2} \langle t_2\rangle$ is isomorphic to the abelianisation of  $G_{ab} \rtimes_{\bar \phi_2} \mathbb{Z} $, therefore  isomorphic to the quotient of  $ (C_1\times \dots \times C_p)$  by the relations $c_i=c_j, i<j<p$.  Since each $c_i$ generates a direct factor of $C_i$, we have the result. 

\end{pf} 

\subsubsection{Structure of the group    $ Out(  G_2 \rtimes \langle t_2\rangle )$}

  A theorem of Guirardel and Levitt  provides a structural feature of    ${\rm Out}(  G_2\rtimes_{\phi_2} \langle t_2\rangle )$, because of the relative hyperbolicity of  $G_2\rtimes_{\phi_2} \langle t_2\rangle $  (from Corollary \ref{coro;toral}).   They first prove (although a modern definitive reference is their essay \cite{VGGL-JSJ})  that there exists a canonical JSJ decomposition for $G_2\rtimes_{\phi_2} \langle t_2\rangle$ as a finite graph of groups in which vertex groups can be surface groups (with boundary), parabolic subgroups, hence free abelian, in which the collection of adjacent edge groups generate a direct factor, and other, so-called rigid groups (see also \cite[\S 10]{DGr}).

\begin{lem}\label{lem;tree_noQH}
  The JSJ graph-of-group decomposition of  $  G_2\rtimes_{\phi_2} \langle t_2\rangle $ contains no vertex group that is a surface groups (with boundary), and has an underlying graph that is a finite tree.
\end{lem}  

\begin{pf}  We refer to  \cite[Prop. 2.11]{DTAMS} for the first claim. For the second claim, take the Bass-Serre tree $T$ of the decomposition. Consider it as a $G_2$-tree, $G$ being a normal subgroup of  $G_2\rtimes_{\phi_2} \langle t_2\rangle$. The quotient $G_2\backslash T$ gives a graph of groups decomposition of $G_2$  that is a free decomposition, i.e. edge groups in $G$ are trivial (see \cite[Lem. 2.8]{DTAMS}). By assumption on $G$ and unicity of the Grushko decomposition, the graph $G\backslash T$ is a tree. The underlying graph of the decomposition of $G\rtimes_{\phi_2} \langle t_2\rangle$ is the quotient of $G\backslash T$ by the action of the cyclic group induced by $\phi_2$. But since $\phi_2$ preserves the conjugacy class of all free factors of $G$, it induces the identity on the graph $G\backslash T$. Thus,  underlying graph of the decomposition of $G\rtimes_{\phi_2} \langle t_2\rangle$ is a tree.
\end{pf}

A graph-of-group presentation for $ G\rtimes_{\phi_2}  \langle t_2\rangle$ is then given by presentations for each vertex and unoriented-edge groups, and   relations given by attaching maps of each oriented edge. %

Guirardel and Levitt  prove in  \cite[Thm 1.4]{VGGL1}   that  there exists a finite index subgroup  ${\rm Out}^1(  G_2\rtimes_{\phi_2} \langle t_2\rangle )$ that preserve the conjugacy class of each maximal parabolic group, and  that fits ain a short exact sequence 

$$  1\to \mathbb{T} \, \longrightarrow \,  {\rm Out}^1(  G_2\rtimes_{\phi_2} \langle t_2\rangle ) \longrightarrow \,     (\prod_{j\in J_S} MCG^0_j)\times ( \prod_{j\in J_P} GL_{n_j,m_j}(\mathbb{Z}))  \to 1   $$

in which $\mathbb{T}$ is a free abelian group,
 and  in which the groups $MCG^0_j$ are Mapping Class Groups for surface vertices of the JSJ decomposition of $ G_2\rtimes_{\phi_2} \langle t_2\rangle  $ (indexed by the set $J_S$), and in which  the groups 
$GL_{n_j,m_j}(\mathbb{Z})$  are the groups of automorphisms of $\mathbb{Z}^{n_j+m_j}$ fixing the first $n_j$ generators.   

As we recalled in the Lemma \ref{lem;tree_noQH}, $J_S$ is empty in our case. 

The factors $GL_{n_j,m_j}(\mathbb{Z})$  correspond to  the automorphism groups of parabolic vertex groups $P_j$  (indexed by the set $J_P$), of rank $n_j+m_j$, fixing the direct factor $E_j<P_j$ generated by  adjacent edge groups, of rank $n_j$. 

The group $\mathbb{T}$ is generated by Dehn twists over edges of the JSJ graph of groups. We will only be interested by its image in the automorphism group of the abelianisation of $ G_2\rtimes_{\phi_2} \langle t_2\rangle  $, and this image is  trivial, since in our setting all Dehn twists are piecewise conjugations (all edges of the graph of groups are separating).

Using  \cite[Thm 4.4]{DGr} on the relatively hyperbolic structures of the vertex groups of the JSJ decomposition,   we can algorithmically compute a collection of coset representatives of  ${\rm Out}^1(  G_2\rtimes_{\phi_2} \langle t_2\rangle )$ in  ${\rm Out}(  G_2\rtimes_{\phi_2} \langle t_2\rangle )$.  Let $\{ \theta_1, \dots, \theta_c\}$ be such a collection. Let $\alpha_i= \theta_i\circ \alpha$. The following is straightforward, from Lemma \ref{orbit}.

\begin{lem} \label{lem;orbit_i} There exists $\psi$ solving the orbit problem of Lemma \ref{orbit} in     ${\rm Aut}(  G_2\rtimes_{\phi_2} \langle t_2\rangle )$ if and only if there exists $i\leq c$, and $\psi_i \in {\rm Aut}(  G_2\rtimes_{\phi_2} \langle t_2\rangle )$ whose class in ${\rm Out}(  G_2\rtimes_{\phi_2} \langle t_2\rangle )$ is in ${\rm Out}^1(  G_2\rtimes_{\phi_2} \langle t_2\rangle )$, such that  the inclusions of  Lemma \ref{orbit} \ref{inclusion} are satisfied for $\alpha_i$ andf $\psi_i$:
 
$$        \bar G_2 \subset   (\overline{\psi_i(\alpha_i (G_1) )} )  \quad \hbox{ and }  \quad  (\overline{\psi_i(\alpha_i (t_1) )} ) \in   \bar t_2  \bar G_2.$$

  \end{lem}

\subsubsection{The orbit problem in the terms of the short exact sequence}

We first agree on a section of $\prod GL_{n_j,m_j}(\mathbb{Z}) $ in ${\rm Out}^1(   G_2\rtimes_{\phi_2} \langle t_2\rangle )$.  This can be done by choosing a graph-of-group presentation of   $G_2\rtimes_{\phi_2} \langle t_2\rangle $  for the JSJ decomposition, and for which  the parabolic vertex groups of the JSJ decomposition are generated by  a certain  basis (as free abelian group) such that the $n_j$ first vectors are in the group generated by the adjacent edge groups. One then realises the group   $GL_{n_j,m_j}(\mathbb{Z})$, as an automorphism group  of the $j$-th parabolic vertex group, that induces the identity on the adjacent edge groups, and therefore, that embeds as a subgroup $\widetilde{GL}_{n_j,m_j}(\mathbb{Z}) $  of the automorphism group of the graph of groups, hence of    ${\rm Out}^1(   G_2\rtimes_{\phi_2} \langle t_2\rangle )$.

Now with this choice, a collection of elements of $\widetilde{GL}_{n_j,m_j}(\mathbb{Z})$    together with  an element in  $\mathbb{T}$
 defines an element of ${\rm Out}^1(   G_2\rtimes_{\phi_2} \langle t_2\rangle )$, and conversely an element in the later group defines a collection in  $\widetilde{GL}_{n_j,m_j}(\mathbb{Z})$ and an element in  $\mathbb{T}$, and all  elements of ${\rm Out}^1(   G_2\rtimes_{\phi_2} \langle t_2\rangle )$ are thus  obtained.

 Let us consider a generating set $\calS$ for $G$, compatible with the free product: $\calS= \bigsqcup_i \calS_{a,i}$  for generating sets $\calS_{a,i}$  for the $A_i$.      
 One can derive  a generating set $ \bar \calS$  of $\bar G_2 $, 
 as an union  $\bigsqcup_j \bar \calS_{a,j} $, where each   $ \bar \calS_{a,j}$  is in the image of the groups $P_j$.

Note that $\widetilde{GL}_{n_j,m_j}(\mathbb{Z}) $ is the identity on all elements of  $ \calS$   except those in  $  \calS_{a,j}$.    Hence it is the identity on all elements of  $ \bar \calS$ except those in  $ \bar \calS_{a,j}$.

\begin{lem}\label{lem;16}
  Fix $i\leq c$.  There exists $[\psi] \in  {\rm Out}^1(  G_2\rtimes_{\phi_2} \langle t_2\rangle )$ satisfying the  inclusions in   
Lemma \ref{lem;orbit_i}   
 if and only if
  for all $j\in J_P$, there exists    $\rho_j\in \widetilde{GL}_{n_j,m_j}(\mathbb{Z})  $ such that  $ \overline{\rho_j} ( x)  \in \overline { \alpha_i( G_1 \cap Q_{j,i}) }$, for all $x\in   \bar \calS_{a,j}$  and for $Q_{j,i}$ the preimage of $P_j$ by $\alpha_i$, and, if  $t_2$ in the abelianisation is in the image of  $P_j$, such that  $\overline{\rho_j}(\bar t_2) \in  \overline{ \alpha_i( t_1{G_1})}  $.   
\end{lem}

\begin{pf}
If there exists $[\psi] \in {\rm Out}^1(  G_2\rtimes_{\phi_2} \langle t_2\rangle )$ resolving the orbit problem  of Lemma \ref{lem;orbit_i} (for $\alpha_i$), then, one consider $\psi^{-1}$ and its  decomposition in $\mathbb{T}$ and $\widetilde{GL}_{n_j,m_j}(\mathbb{Z})  $.  The automorphism $\psi^{-1}$ sends $\bar G_2$ to the image of $\bar G_1$ by $\alpha_i$.  One easily derives the given relations since $\mathbb{T}$ acts trivially on the abelianisation.

Conversely, it there are $\rho_j$      
as in the statement, since their support among the elements  of $ \bar \calS$  are disjoint, they define an element of  $ {\rm Out}^1(  G_2\rtimes_{\phi_2} \langle t_2\rangle )$ sending    $\bar G_2$ to the image of $\bar G_1$ by $\alpha_i$.  Its inverse thus satisfies  the orbit problem of Lemma \ref{lem;orbit_i}.

\end{pf}

\begin{lem} \label{lem;17}  For each $j\in J_P$, it is decidable whether there exists an element $\rho_j$ of  $\widetilde{GL}_{n_j,m_j}(\mathbb{Z})$ such that, 
 for all $x\in   \bar \calS_{a,j}$, 
 $ \overline{\rho_j}  (x)  \in \overline { \alpha_i( G_1 \cap Q_{j,i}) }$, and,  if  $t_2$ in the abelianisation is in the image of  $P_j$, such that  $\overline{\rho_j}(\bar t_2) \in  \overline{ \alpha_i( t_1{G_1})}  $.    
\end{lem}

\begin{pf} 
This amounts to check whether there is a matrix in $GL_{n_j,m_j}(\mathbb{Z})$ sending an hyperplane of $\mathbb{Z}^{n_j+m_j}$  in another, and possibly, an element in a given hyperplane, thus on a point in the suitable quotient.  
The first condition,  expressed in the dual, is requiring sending a vector on another.

 The matrices   in $GL_{n_j,m_j}(\mathbb{Z})$ are exactly those that are with integer entries, triangular by block, with the top diagonal block being the identity, the bottom  diagonal block having determinant $\pm 1$, and the top right block arbitrary. This is an arithmetic group in an explicit algebraic subgroup of matrices.  One can thus  apply  Grunewald and Segal's algorithm,  \cite[Algorithm A]{GS} in order to treat this orbit problem.

\end{pf}

By Lemma \ref{lem;17} and \ref{lem;16}, the orbit problem of Lemma \ref{lem;orbit_i} is decidable. Therefore, the problem of Lemma \ref{orbit} is decidable. Therefore,  the conjugacy problem in $Out(G,\calG)$  for atoroidal toral automorphisms of $G$ is then solved.


{\small 

{\footnotesize
{\sc \noindent Francois Dahmani, Institut Fourier, Univ. Grenoble Alpes, F-38000 Grenoble \\
    {\tt francois.dahmani@univ-grenoble-alpes.fr}\\ 
    \noindent Ruoyu Li,  School of Mathematics, Jiaying University, Meijiang District, Meizhou, Guangzhou, China. \\
{\tt 201901113@jyu.edu.cn}
} }
}

\end{document}